\documentclass[fleqn,reqno,11pt,a4paper,final]{amsart}

\usepackage[a4paper,left=20mm,right=20mm,top=20mm,bottom=20mm,marginpar=20mm]{geometry}
\usepackage{amsmath}
\usepackage{amssymb}
\usepackage{amsthm}
\usepackage{amscd}
\usepackage[ansinew]{inputenc}
\usepackage[noadjust]{cite}
\usepackage{bbm}
\usepackage{xcolor}
\usepackage[english=american]{csquotes}
\usepackage[final]{graphicx}
\usepackage{hyperref}
\usepackage{calc}
\usepackage{mathptmx}
\usepackage{mathtools}
\usepackage{tikz}
\usepackage{graphicx}
\usepackage{stix}
\usepackage{soul}
\usepackage{ulem}
\usepackage{hyperref}
\hypersetup{
    colorlinks,
    linkcolor={luh-dark-blue},
    citecolor={luh-dark-blue},
    urlcolor={luh-dark-blue}
}

\graphicspath{{../Pictures/}}

\numberwithin{equation}{section}

\newtheoremstyle{thmlemcorr}{10pt}{10pt}{\itshape}{}{\bfseries}{.}{10pt}{{\thmname{#1}\thmnumber{ #2}\thmnote{ (#3)}}}
\newtheoremstyle{thmlemcorr*}{10pt}{10pt}{\itshape}{}{\bfseries}{.}\newline{{\thmname{#1}\thmnumber{ #2}\thmnote{ (#3)}}}
\newtheoremstyle{remexample}{10pt}{10pt}{}{}{\bfseries}{.}{10pt}{{\thmname{#1}\thmnumber{ #2}\thmnote{ (#3)}}}
\newtheoremstyle{ass}{10pt}{10pt}{}{}{\bfseries}{.}{10pt}{{\thmname{#1}\thmnumber{ A#2}\thmnote{ (#3)}}}

\theoremstyle{thmlemcorr}
\newtheorem{theorem}{Theorem}
\numberwithin{theorem}{section}
\newtheorem{lemma}[theorem]{Lemma}
\newtheorem{corollary}[theorem]{Corollary}

\theoremstyle{thmlemcorr*}
\newtheorem{theorem*}{Theorem}
\newtheorem{lemma*}[theorem]{Lemma}
\newtheorem{corollary*}[theorem]{Corollary}
\newtheorem{proposition*}[theorem]{Proposition}
\newtheorem{problem*}[theorem]{Problem}
\newtheorem{conjecture*}[theorem]{Conjecture}
\newtheorem{definition*}[theorem]{Definition}
\newtheorem{assumption*}[theorem]{Assumption}

\theoremstyle{remexample}
\newtheorem{remark}[theorem]{Remark}

\theoremstyle{ass}

\newcommand{\Acal}{\mathcal{A}}

\newcommand{\Fcal}{\mathcal{F}}

\newcommand{\Hcal}{\mathcal{H}}

\newcommand{\Lcal}{\mathcal{L}}

\newcommand{\longhookrightarrow}{\lhook\joinrel\longrightarrow}

\newcommand{\ii}{\mathrm{i}}

\newcommand{\N}{\mathbb{N}}
\newcommand{\R}{\mathbb{R}}
\newcommand{\C}{\mathbb{C}}

\newcommand{\eps}{\epsilon}

\definecolor{luh-dark-blue}{rgb}{0.0, 0.313, 0.608}

\def\XXint#1#2#3{{\setbox0=\hbox{$#1{#2#3}{\int}$}
\vcenter{\hbox{$#2#3$}}\kern-.5\wd0}}

\renewcommand{\eps}{\varepsilon}
\renewcommand{\epsilon}{\varepsilon}
\renewcommand{\phi}{\varphi}

\begin{document}

%% TITLE MATTERS

\title[]{Non-Newtonian two-phase thin-film problem: Local existence, uniqueness, and stability} % with Ellis law}

\author{Oliver Assenmacher}
\address{\textit{Oliver Assenmacher:}  Institute of Applied Mathematics, University of Bonn, Endenicher Allee~60, 30167 Bonn, Germany}
\email{s6olasse@uni-bonn.de}

\author{Gabriele Bruell}
\address{\textit{Gabriele Bruell:}Institute for Analysis, Karlsruher Institute of Technology (KIT), D-76128 Karlsruhe, Germany}
\email{gabriele.bruell@kit.edu}

\author{Christina Lienstromberg}
\address{\textit{Christina Lienstromberg:}  Institute of Applied Mathematics, University of Bonn, Endenicher Allee~60, 30167 Bonn, Germany}
\email{lienstromberg@ifam.uni-hannover.de}

\begin{abstract}
We study the flow of two immiscible fluids located on a solid bottom, where the lower fluid is Newtonian and the upper fluid is a non-Newtonian Ellis fluid. Neglecting gravitational effects, we consider the formal asymptotic limit of small film heights in the two-phase Navier--Stokes system. This leads to a strongly coupled system of two parabolic equations of fourth order with merely H\"older-continuous dependence on the coefficients. For the case of  strictly positive initial film heights we prove local existence of strong solutions by abstract semigroup theory. Uniqueness is proved by energy methods. 
% {\color{orange} \sout{Finally, we use center manifold theory to show that solutions exist forever and converge exponentially fast to a flat steady state, provided they are close to the steady state initially.  }
Under additional regularity assumptions, we  investigate asymptotic stability of the unique equilibrium solution, which is given by constant film heights.
%We study the problem of existence and uniqueness of strong solutions to a two-phase thin film problem with one non-Newtonian Ellis fluid on top of a Newtonian fluid. 
%
%Under the assumption that capillarity is the only driving force, we consider the formal asymptotic limit of small film heights in the two-phase Navier--Stokes system. This leads to a strongly coupled system of two parabolic equations of fourth order with only H\"older-continuous dependence on the coefficients. For the non-degenerate case of positive film heights we prove local existence of strong solutions by abstract semigroup theory. Uniqueness is proved by energy methods. Finally, we prove that solutions exist forever and converge exponentially fast to a flat steady state, provided they are close to the steady state initially.  
\end{abstract}
\vspace{4pt}

%\noindent\textsc{MSC (2010): XXXXX (primary);}

%\noindent\textsc{Keywords:}

%\vspace{4pt}

%\noindent\textsc{Date:} \today{}.
%\end{abstract}

%% PDF MATTERS

%% START OF CONTENT

\maketitle

\section{Introduction}

\subsection*{Aim of the paper}
The present manuscript is devoted to the mathematical modeling and analysis of a two-phase stratified thin-film problem, where one of the fluids features a shear rate-dependent viscosity. 
More precisely, we assume that the two-phase flow is located on an impermeable bottom and that the lower fluid is Newtonian, while the upper fluid is non-Newtonian with a viscosity $\mu$ determined by the Ellis constitutive law \cite{WS:1994}. For an Ellis fluid the viscosity is implicitly given by the relation
\begin{equation}\label{eq:Ellis}
	\frac{1}{\mu} = \frac{1}{\mu_0} \left(1 + \left|\frac{\tau}{\tau_{1/2}}\right|^{p-1}\right), \quad p \geq 1,
\end{equation}
where $\tau$ is the shear stress, $\mu_0$ denotes the viscosity at zero shear stress and $\tau_{1/2}$ is the shear stress at which the viscosity has dropped to $\mu_0/2$. Assuming that $\tau_{1/2} \in (0,\infty)$, the choice of the flow-behaviour exponent $p > 1$ reflects the shear-thinning behaviour of the fluid. Moreover, one can recover a Newtonian behaviour for $p=1$ or $1/\tau_{1/2}\to 0$. The fact that the Ellis law is able to describe a shear-thinning behavior for moderate shear stresses on the one hand and a Newtonian plateau for very low shear stresses on the other hand makes it a more realistic model for many real fluids, compared to the %well-known
widely-used power-law \cite{Myres} for instance. 
%{\color{orange} GB: I think the Ellis law does only reflect the plateau behavior for low sheer stress, not for large. A model describing both plateaus is for instance given by the Carreau model, cf. \cite{Myres}\footnote{In \cite[page 2 (article in Dropbox)]{Myres}: \textit{This model [Ellis] cannot predict the second Newtonian plateau, however, as stated this is usually not of interest and particularly not for this investigation where it is the low shear region that is of primary interest.}}.} See Figure \ref{F:1} for the geometrical setting.

\begin{center}
\begin{figure}[h]\label{fig1}
\begin{tikzpicture}[domain=0:3*pi, xscale=1.2, yscale=0.8] 
\draw[ultra thick, smooth, variable=\x, luh-dark-blue!70] plot (\x,{0.3*cos(\x r)+2}); 
\fill[luh-dark-blue!10] plot[domain=0:3*pi] (\x,{0.2*sin(\x r)+1}) -- plot[domain=3*pi:0] (\x,{0.3*cos(\x r)+2});
\draw[ultra thick, smooth, color=luh-dark-blue] plot (\x,{0.2*sin(\x r)+1});
\fill[luh-dark-blue!20] plot[domain=0:3*pi] (\x,0) -- plot[domain=3*pi:0] (\x,{0.2*sin(\x r)+1});
\draw[very thick,<->] (3*pi+0.4,0) node[right] {$x$} -- (0,0) -- (0,3.5) node[above] {$z$};
%\node[right] at (3*pi,1.7) {$g(t,x)$};
%\node[right] at (3*pi,0.9) {$f(t,x)$};
%  \draw (2*pi,1.5) node {$\Omega_+(t)$}; 
%     \draw (2*pi,0.4) node {$\Omega_-(t)$}; 
\draw[-] (0,-0.3) -- (0.3, 0);
\draw[-] (0.5,-0.3) -- +(0.3, 0.3);
\draw[-] (1,-0.3) -- +(0.3, 0.3);
\draw[-] (1.5,-0.3) -- +(0.3, 0.3);
\draw[-] (2,-0.3) -- +(0.3, 0.3); 
\draw[-] (2.5,-0.3) -- +(0.3, 0.3);
\draw[-] (3,-0.3) -- +(0.3, 0.3);
\draw[-] (3.5,-0.3) -- +(0.3, 0.3);
\draw[-] (4,-0.3) -- +(0.3, 0.3);
\draw[-] (4.5,-0.3) -- +(0.3, 0.3);
\draw[-] (5,-0.3) -- +(0.3, 0.3);
\draw[-] (5.5,-0.3) -- +(0.3, 0.3);
\draw[-] (6,-0.3) -- +(0.3, 0.3);
\draw[-] (6.5,-0.3) -- +(0.3, 0.3);
\draw[-] (7,-0.3) -- +(0.3, 0.3);
\draw[-] (7.5,-0.3) -- +(0.3, 0.3);
\draw[-] (8,-0.3) -- +(0.3, 0.3);
\draw[-] (8.5,-0.3) -- +(0.3, 0.3);
\draw[-] (9,-0.3) -- +(0.3, 0.3);
\end{tikzpicture}   
\caption{Two-phase flow on an impermeable bottom. The lower fluid is Newtonian, while the upper one is non-Newtonian.}\label{F:1}
\end{figure}
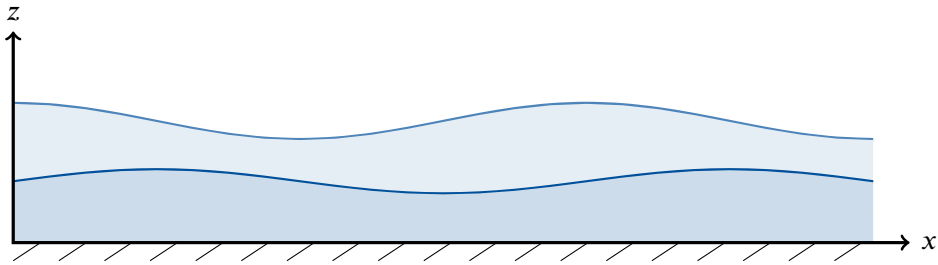 
\end{center}
{Besides the non-Newtonian rheological behaviour of the upper fluid, two assumptions are crucial in the derivation of the model. First, we assume the flow to be uniform in one horizontal direction such that it depends on one spatial variable only. Second, we assume that the two liquid films are very thin and consider the asymptotic limit of vanishing heights. Therefore, starting with a Navier--Stokes system for the two immiscible fluids, we apply lubrication approximation and cross-sectional averaging in order to derive a coupled system of evolution equations for the two film heights $f=f(t,x)$ and $g=g(t,x)$ of the lower Newtonian and the upper non-Newtonian fluid, respectively.} 
%Assuming the flow to be uniform in one horizontal direction, we apply lubrication approximation and cross-sectional averaging to the governing equations 
%in order to derive a system of coupled order evolution equations for the two film heights $f=f(t,x)$ and $g=g(t,x)$ of the lower Newtonian and the upper non-Newtonian fluid, respectively. 
The functions $f,g$ depend on the temporal variable $t>0$ and the spatial variable $x\in \Omega$, where $\Omega\subset \R$ is a finite interval. 
%If capillarity is the main driving force, the evolution equations are of fourth order and given by
{
Neglecting gravitational effects and assuming that capillarity is the system's driving force, we obtain the coupled system}
\begin{flalign}\label{eq:system_I}
 	\begin{cases}
 		f_t + \left(m s^+ \left( \frac{f^3}{3}+ \frac{f^2 g}{2} \right) (f+g)_{xxx}
 		+\frac{s^-}{3} f^3 f_{xxx}\right)_x = 0, 
 		 & 
 		\\
 		g_t + 
 		\left(\mkern-6mu\left(\mkern-3mu\frac{ms^+}{2} f^2g + ms^+ f g^2 + \frac{s^+}{3 \mu_0^+} g^3 \mkern-3mu\right)(f+g)_{xxx} + \frac{s^-}{2} f^2 g f_{xxx} + C_p |g|^{p + 2} \Phi\bigl((f+g)_{xxx}\bigr)\mkern-6mu\right)_{\mkern-6mu x}
 		\mkern-3mu = 0, & 
 	\end{cases}&&
\end{flalign}
{of fourth order that describes the dynamics of the two interfaces $f$ and $g$. Here, $s^\pm$ is the constant dimensionless surface tension of the lower (-) and the upper (+) fluid, respectively. Moreover, $\mu_0^\pm > 0$ is the characteristic viscosity of the lower and the upper fluid at zero shear, respectively, and $m = \mu_0^+/\mu_0^- > 0$ is the ratio of the two characteristic viscosities. The function $\Phi$ in \eqref{eq:system_I} is related to the non-Newtonian part in the Ellis law and is defined by
\begin{equation*}
	\Phi(d) = |d|^{p-1}d, \quad d \in \R,
\end{equation*} 
where $p>1$ is the flow behaviour exponent and $C_p > 0$ is a positive constant that depends on $p$.
%where $m,s^\pm$ and $\mu_0>0$ are constants related to the zero stress viscosity and surface tension. Moreover, $p\geq 1$, the constant $C_p>0$ and the function $\Phi(d)=|d|^{p-1}d$ are related to the imposed Ellis-law for the upper fluid.

For the coupled system \eqref{eq:system_I}, supplemented with the Neumann-type boundary conditions
\[
 f_x = f_{xxx}=g_x = g_{xxx}=0\quad \mbox{at } \partial \Omega,
\]
and initial conditions
\[
f(0,x)=f_0(x),\quad g(0,x)=g_0(x),\quad x\in \Omega,
\]
we prove local in time existence of strong solutions, which are unique with respect to the initial data. 
%Note that due to the degeneracy of the system when the film heights $f,g$ drop to zero, the existence of global strong solutions is in general out of reach. 
We obtain the existence of local positive strong solutions for strictly positive initial data  by %parabolic 
analytic semigroup theory. Since the coefficient function $\Phi$ is merely $(p-1)$-H\"older continuous for $p\in (1,2)$, uniqueness of solutions with respect to initial data does not follow by the classical abstract theory of {\sc{Amann}} \cite{A:1993}. Instead, we use energy methods to prove a uniqueness result for local strong solutions.

Note that due to the degeneracy of the system when the film heights $f,g$ drop to zero, the existence of global strong solutions can in general not be expected.
%is in general out of reach. 
However, we use center manifold theory \cite{HI,Milke} in order to investigate the long-time behaviour of solutions which are initially close to constant positive film heights. 
%Furthermore, we use center manifold theory \cite{HI,Milke} in order to investigate the asymptotic stability of steady states. 
In view of the boundary conditions one may observe that flat films of positive height are the only possible steady state solutions of \eqref{eq:system_I}. We prove that, if the film heights are initially close to positive flat films, then the corresponding solution exists globally in time and converges  exponentially to the two-phase flat film as time tends to infinity. In other words, for large times, surface tension effects dominate the characteristic stresses of the Ellis fluid and solutions behave as in the Newtonian regime.

\subsection*{Non-Newtonian fluids}
%According to the concrete application different types of fluids are considered
Although many common liquids and gases, such as water and air, may reasonably be considered Newtonian,  there is still a multitude of real fluids which are in fact non-Newtonian.
%considered in real situations are in fact non-Newtonian}.
Newtonian fluids obey Newtons law of viscosity, that is, they are characterised by a linear dependence of the shear stresses on the local strain rate. This dependence is nonlinear for {non-Newtonian} fluids. In other words, %the viscosity of 
non-Newtonian fluids become more solid or more liquid under force. Among others, two prominent classes of shear rate-dependent fluids are {shear-thinning} (pseudoplastic) and {shear-thickening} (dilatant) fluids. Shear-thinning fluids are characterised by a viscosity that decreases with increasing shear rate, as for instance wall paint. On the contrary, shear-thickening fluids exhibit a viscosity that increases with increasing shear rate, as for instance %cornflower 
corn starch mixed with water. 
%Fluids with a purely shear-thickening or purely shear-thinning behaviour are called Ostwald--de Waele power-law fluids. 
A famous class of fluids with a purely shear-thickening or purely shear-thinning behaviour consists of so-called Ostwald--de Waele fluids or power-law fluids.
%For this type of fluid the shear-stress $\tau$ is given by $\tau(s) = \mu_0 \mu(s)s$, with an effective viscosity $\mu(s)=|s|^{q-1}$, a characteristic viscosity $\mu_0 \in \R_+$, and a flow behaviour exponent $q \in \R_+$.
For this type of fluid the viscosity $\mu$ is given by 
\begin{equation} \label{eq:power-law_intro}
	\mu = \mu_0 \left|\tau\right|^\frac{q-1}{q} 
\end{equation}
with a characteristic viscosity $\mu_0 \in \R_+$ and a flow behaviour exponent $q \in \R_+$.
Though convenient to use because of its simple mathematical structure, the validity of a power-law \eqref{eq:power-law_intro} is  often restricted to certain ranges of shear rates. In other words,  the power-law model might lead to inaccurate predictions at very high or very low shear rates. In the shear-thinning regime $q<1$, for instance, an Ostwald--de Waele fluid would have zero viscosity as the shear rate tends to infinity and infinite viscosity at rest.

Many fluids display a combination of different rheological behaviours at different ranges of shear rates. Among polymeric fluids, for instance,  there is a variety that has a shear-thinning behaviour for moderate shear rates but a Newtonian behavior for very high and very low shear rates, respectively. More precisely, at low and high shear rates, respectively, the viscosity is almost constant, while at moderate shear rates the rheological behaviour can be modeled by a power-law. 

Regarding this phenomenon, a more  promising model than the power-law is the Ellis constitutive law, {cf. \eqref{eq:Ellis}}, which captures a Newtonian behavior close to zero shear and a shear-thinning power-law behaviour for higher shear rates \cite{Myres,WS:1994}. We refer to Figure \ref{F:compare} for a schematic comparison of the viscous behaviour of a Newtonian fluid, a power-law fluid and an Ellis fluid depending on the shear rate.

\begin{center}
\begin{figure}[h!]
\begin{tikzpicture}[yscale=1.5, xscale=1.5]
\draw[->, very thick] (0,0)--(4,0);
\draw[->, very thick] (0,0)-- (0,2.5);
\node at (2,-0.2) {shear rate};
\node[rotate=90] at (-0.2,1.2) {viscosity};
\draw[-, dashed, gray, thick] (0,1.5)--(4,1.5) node[above left]{\color{black}{Newtonian}};
\draw[domain=0.2:4,variable=\x,red!40, very thick] plot ({\x},{(1.5-0.3*abs(\x))^2 });
\draw[domain=0:3,variable=\x,luh-dark-blue, very thick] plot ({\x},{1.5-0.13*(\x+0.3) * (\x+0.3) });
\node at (3.5,0.5) {power-law};
\node at (0.6,1.1) {Ellis-law};
\end{tikzpicture}
\caption{ Common models for the rheological behavior of shear-thinning fluids: Power-law (red) and Ellis-law (blue). The dashed line represents the Newtonian law, where the viscosity does not depend on the shear rate. }\label{F:compare}
\end{figure}
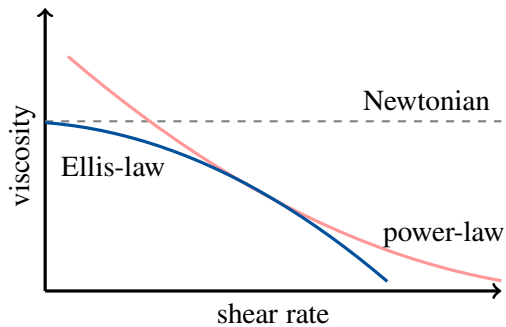
\end{center}

\subsection*{Two-phase thin films}
%A stratified two-phase or even multi layered thin film consists of two or more layers of  viscous fluids on a solid bottom with thickness typically ranging from a nanometer to several micrometers. 
A stratified two-phase thin film consists of two layers of immiscible viscous fluids on a solid bottom with thicknesses typically ranging from a nanometer to several micrometers.}
They find applications in various areas, such as medical sciences, chemistry, or industrial processes. 
Examples include techniques for concentrating leukocytes from blood \cite{SW} in medical treatments or  industrial multi layer coating possess \cite{RW}. 
An example for a stratified two-phase thin film consisting of a non-Newtonian fluid on top of a Newtonian one, as studied in the present manuscript, is given by assembling artificial membranes. Artificial membranes are thin layers of commonly non-Newtonian (often shear-thinning) fluids \cite{Espinosa,Hermans} mimicking cell membranes in nature and find applications in nanotechnology and biosensors \cite{CC}.
% Furthermore, one commonly used method of assembling artificial membranes, so called solid supported lipid bilayers, is to deposit them on a solid substrate in an aqueous environment, which leads to the formation of a thin water layer beneath the membrane. 
The main advantage of considering an artificial membrane on a thin film of water, which is considered as a Newtonian fluid, compared to freely-suspended membranes, is an increase in stability, allowing for more avenues of experimentation \cite{Siontorou}.

Commonly, the mathematical model for the evolution of a single or multi layered thin film flow is derived from the Navier--Stokes equations in each constituent fluid region by so-called \textit{lubrication approximation} and cross-sectional averaging. 
%Depending on the geometrical setting, a boundary condition on the impermeable bottom and stress balance conditions at the fluid--fluid and fluid--air interfaces have to be imposed.
Depending on the particular geometrical setting, %\textcolor{magenta}{particular physical assumptions}, 
boundary conditions on the impermeable bottom and stress balance conditions at the fluid--fluid and fluid--air interfaces have to be imposed.
In thin film systems effects due to viscosity, diffusion, and surface tension become ever more important and may dominate gravitational or inertial forces.

%{\color{gray}One application of thin-films are artificial membranes. These are thin bilayers of lipids mimicking the com- position of cell membranes in nature. They are studied to help understand the behaviour and function of real cell membranes but are also of interest due to applications in biosensors and nanotechnology. Depending on the external conditions, i.e. pressure and temperature, these lipid layers behave similar to viscous liquids and some lipids have been observed to exhibit shear-thinning behaviour.
%Furthermore, one commonly used method of assembling artificial membranes, so called solid supported lipid bilayers, is to deposit them on a solid substrate in an aqueous environment, which leads to the formation of a thin water layer beneath the membrane. The main advantage of this procedure compared to freely-suspended membranes, which consist of a lipid bilayer spanning a small hole between two containers filled with an aqueous solution, is an increase in stability, allowing for more avenues of experimentation.}

\subsection*{Related results.} 
%Commonly, the mathematical model for the evolution equations of thin films is derived from the Navier--Stokes equations by lubrication approximation and cross-sectional averaging.
Assuming the fluid to be located on an impermeable flat bottom and uniform in one horizontal direction, 
the classical thin-film equation for Newtonian fluids reads
\begin{equation}\label{eq:thin_film}
	h_t+\left(h^nh_{xxx}\right)_x=0,\qquad n\in \N,
\end{equation}
where $h=h(t,x)$ represents the thickness of the thin film, depending on the temporal variable $t>0$ and spatial variable $x \in \Omega$, where $\Omega \subset \R$ is a finite interval. The case $n=1$ corresponds to %the flow of a Hele--Shaw cell 
a Hele-Shaw flow, while  $n=2,3$ model the flow of a capillary driven Newtonian fluid with either Navier--slip boundary condition ($n=2$) or no-slip boundary condition ($n=3$) imposed on the solid bottom.
%flow \textit{Navier--slip} and  \textit{no-slip}  condition on the bottom, respectively. 
{There exists a rich literature on the dynamics of a single Newtonian thin film as in \eqref{eq:thin_film}.}
Pioneering works on the existence of global weak solutions of the classical thin-film equations \eqref{eq:thin_film} and their properties are due to {\sc{Bernis}} and {\sc{Friedman}} \cite{BF} , followed by \cite{BBD} of {\sc{Beretta}} et al., and \cite{BP} by {\sc{Bertozzi}} and {\sc{Pugh}}.
%due to {\sc{Bernis}} and {\sc{Friedman}} \cite{BF} followed by {\sc{Beretta}} et al. \cite{BBD} and {\sc{Bertozzi}} and {\sc{Pugh}} \cite{BP}. 
While the paper \cite{BF} covers also uniqueness of weak solutions for $n\geq 4$, it is still one of the probably most famous open problems in the context of \eqref{eq:thin_film} to prove uniqueness of weak solutions for $n<4$. Another issue that has hitherto attracted much attention is to control the solution at the contact points of liquid, solid and gas. Since weak solutions have in general only weak regularity properties, a lot of research has also been dedicated to the verification of existence and uniqueness of non-negative strong solutions of the free-boundary value problem associated to \eqref{eq:thin_film}. We refer the reader to the works \cite{GGKO:2014,GK:2010,GKO:2008,GnPet:2018,Kn:2011,Kn:2015} for strong solutions in (weighted) Sobolev or H\"older spaces with different prescribed slip conditions.
Concerning stratified two-phase generalisations of the thin film equation \eqref{eq:thin_film} for Newtonian fluids, we refer to \cite{BG, EMM, EM}.

In the setting of non-Newtonian one-phase fluids obeying a power-law, we refer to the early works of {\sc{King}} \cite{K1,K2}, where non-Newtonian generalizations of \eqref{eq:thin_film} of the form
\begin{equation}\label{eq:thin_film_nN}
  h_t+\left(h^n|h_{xxx}|^{p-2}h_{xxx}\right)=0,\qquad p> 2,
\end{equation}
are considered. For $p>2$ and $\frac{p-1}{2}<n<2p-1$ {\sc{Ansini}} and {\sc{Giacomelli}} \cite{AG:2004} establish the existence of non-negative global weak solutions and proved the existence of travelling waves of \eqref{eq:thin_film_nN}. 
Considering 
%instead of a power-law 
an Ellis fluid on an impermeable bottom with no-slip boundary condition, the thin-film approximation is given by
\begin{equation}\label{eq:thin_film_Ellis}
h_t +\left( h^3(1+|hh_{xxx}|^{p-1})h_{xxx} \right)_x=0,\qquad p\geq 1.
\end{equation}
If $p=1$, the above equation reduces to the classical thin-film equation for Newtonian fluids \eqref{eq:thin_film} for $n=3$. In \cite{LM:2020} the authors prove the existence of local strong solutions of \eqref{eq:thin_film_Ellis} and uniqueness with respect to initial data. Moreover, the authors of \cite{AG:2002} investigate a class of quasi-self-similar solutions describing the spreading of a droplet in the limit of a Newtonian rheology. Exploiting the non-Newtonian behaviour only in a small region near the contact line, the authors observe that the usage of a shear-thinning rheology avoids the well-known no-slip paradox.

Note that the two-phase system \eqref{eq:system_I} reduces formally to the one-phase non-Newtonian thin-film equation \eqref{eq:thin_film_Ellis} if we set $f=0$.

\subsection*{Outline of the paper.} We close the introduction with a brief outline of the paper. 
In Section \ref{sec:modelling} we derive the evolution equations for the fluid-fluid interface and the fluid-air interface, respectively. This is done by taking the formal asymptotic limit of vanishing film heights in the Navier-Stokes system for a two-phase problem with general dynamic viscosity functions. Furthermore, we state explicitly the equations for the case in which a non-Newtonian Ellis fluid is placed on top of a Newtonian fluid. Section \ref{sec:existence} is devoted to the proof of the existence of strong solutions for possibly short times. Uniqueness is treated separately in Section \ref{sec:uniqueness}. Finally, in Section \ref{sec:long-time} we restrict to Ellis fluids with flow-behaviour exponents $p\geq 2$ and study the long-time behaviour of solutions which are initially close a flat film. We prove that these solutions have an infinite lifetime and converge to flat films of positive heights as time tends to infinity. 

\bigskip

\section{Modelling}\label{sec:modelling}

\subsection{The Navier--Stokes system}

We study the evolution of two incompressible liquid films on top of each other when capillarity is the only driving force. In particular, gravitational effects are ignored. {Let us recall the geometrical setting.}
The lower fluid is located on an flat solid impermeable bottom at height $z=0$. Moreover, it is assumed to be Newtonian, i.e. it has constant viscosity $\mu^{-} > 0$. On top of this fluid film, we consider a second film that is assumed to be non-Newtonian with a shear-rate dependent viscosity $\mu^+$.
%Its shear-thinning rheology is described implicitly by the so-called Ellis constitutive law \textcolor{magenta}{[Ref]}
%\begin{equation*}
%	\frac{1}{\mu^{+}(|u_z^+|)} = \frac{1}{\mu_0^+}\left(1 + \left|\frac{\mu^+(|u_z^+|) u_z^+}{\tau_{1/2}^+}\right|^{\beta-1}\right).
%	\quad
%%	\tau^2 = \mu^2(|u_z^2|) u_z^2
%\end{equation*}
%Here $\beta > 1$ is the flow-behaviour exponent, $\mu_0^+ > 0$ is the viscosity of the upper fluid at zero shear, and $\tau_{1/2}^+ > 0$ is the shear stress at which the viscosity $\mu^{+}$ has dropped to $\mu_0^+/2$. 
%Moreover, $u^+$ is the horizontal component of the velocity field $\textbf{u}^+$ of the upper fluid and $u^+_z$ is its partial derivative with respect to the vertical variable $z$.
The two fluids are immiscible with a separated interface given by the graph $z=h^-(t,x) > 0$, where $t \geq 0$ denotes the time variable and $x \in \Omega \subset \R$ the horizontal spatial variable. Here, $\Omega \subset \R$ is a bounded interval. The upper fluid is separated from air by the interface $z=h^+(t,x)>h(t,x)^-$. We use the notation
\begin{equation*}
	\begin{cases}
		\Omega_{-}(t) = \left\{(x,z) \in \R^2;\ x\in\Omega,\ 0 < z < h^-(t,x)\right\}, &
		\\
		\Omega_{+}(t) = \left\{(x,z) \in \R^2;\ x\in\Omega,\ h^-(t,x) < z < h^+(t,x)\right\}, &
	\end{cases}
\end{equation*}
to describe the domains occupied by the lower, respectively upper fluid, at any instant $t \geq 0$ of time, cf. Figure~\ref{F:3}.
%\vspace{1cm}
\begin{center}
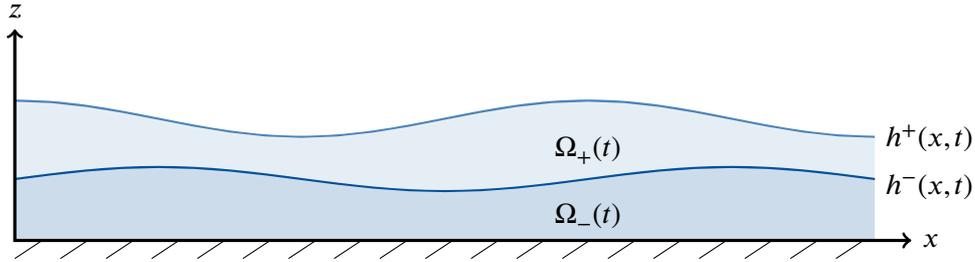
\begin{figure}[h]\label{fig1}
\begin{tikzpicture}[domain=0:3*pi, xscale=1.2, yscale=0.8] 
\draw[ultra thick, smooth, variable=\x, luh-dark-blue!70] plot (\x,{0.3*cos(\x r)+2}); 
\fill[luh-dark-blue!10] plot[domain=0:3*pi] (\x,{0.2*sin(\x r)+1}) -- plot[domain=3*pi:0] (\x,{0.3*cos(\x r)+2});
\draw[ultra thick, smooth, color=luh-dark-blue] plot (\x,{0.2*sin(\x r)+1});
\fill[luh-dark-blue!20] plot[domain=0:3*pi] (\x,0) -- plot[domain=3*pi:0] (\x,{0.2*sin(\x r)+1});
\draw[very thick,<->] (3*pi+0.4,0) node[right] {$x$} -- (0,0) -- (0,3.5) node[above] {$z$};
\node[right] at (3*pi,1.7) {$h^+(x,t)$};
\node[right] at (3*pi,0.9) {$h^-(x,t)$};

 \draw (2*pi,1.5) node {$\Omega_+(t)$}; 
    \draw (2*pi,0.4) node {$\Omega_-(t)$}; 
\draw[-] (0,-0.3) -- (0.3, 0);
\draw[-] (0.5,-0.3) -- +(0.3, 0.3);
\draw[-] (1,-0.3) -- +(0.3, 0.3);
\draw[-] (1.5,-0.3) -- +(0.3, 0.3);
\draw[-] (2,-0.3) -- +(0.3, 0.3);
\draw[-] (2.5,-0.3) -- +(0.3, 0.3);
\draw[-] (3,-0.3) -- +(0.3, 0.3);
\draw[-] (3.5,-0.3) -- +(0.3, 0.3);
\draw[-] (4,-0.3) -- +(0.3, 0.3);
\draw[-] (4.5,-0.3) -- +(0.3, 0.3);
\draw[-] (5,-0.3) -- +(0.3, 0.3);
\draw[-] (5.5,-0.3) -- +(0.3, 0.3);
\draw[-] (6,-0.3) -- +(0.3, 0.3);
\draw[-] (6.5,-0.3) -- +(0.3, 0.3);
\draw[-] (7,-0.3) -- +(0.3, 0.3);
\draw[-] (7.5,-0.3) -- +(0.3, 0.3);
\draw[-] (8,-0.3) -- +(0.3, 0.3);
\draw[-] (8.5,-0.3) -- +(0.3, 0.3);
\draw[-] (9,-0.3) -- +(0.3, 0.3);
\end{tikzpicture}   
\caption{The fluid-air interface $h^+$ and the fluid-fluid interface $h^-$.}\label{F:3}
\end{figure}
\end{center}
%\vspace{1cm}
The dynamics of the two fluids in $\Omega_\pm(t)$  may be described by the (non-Newtonian) Navier--Stokes system 
\begin{equation*}
\begin{cases}
	\rho^{\pm} \bigl(\textbf{u}_t^{\pm} + (\textbf{u}^{\pm}\cdot \nabla) \textbf{u}^{\pm}\bigr)
	=
	\nabla\cdot \textbf{T}^\pm 
	& 
	\text{in } \Omega_{\pm}(t)  \\
	\nabla\cdot \textbf{u}^\pm
	=
	0
	& 
	\text{in } \Omega_{\pm}(t),
\end{cases}
\end{equation*}
where $\textbf{u}^\pm(t,x,z) = (u^\pm(t,x,z),v^\pm(t,x,z))$ 
denotes the velocity field, $p^\pm = p^\pm(t,x,z)$ the pressure, and $\rho^\pm > 0$ the density of each of the fluids. Furthermore, $\textbf{T}^\pm = -p^\pm I + 2\mu^\pm\left(\left\|\textbf{D}^\pm\right\|^2\right) \textbf{D}^\pm$ and $\mu^\pm$ are the stress tensor and the dynamic viscosity of the respective fluid, with $\textbf{D}^\pm = \frac{1}{2}\left(\nabla \textbf{u}^\pm + (\nabla \textbf{u}^\pm)^T\right)$ denoting the symmetric gradient of $\textbf{u}^\pm$ and $\left\|\textbf{D}^\pm\right\| = \sqrt{2 \text{tr}(|\textbf{D}^\pm|^2)}$.
This system is complemented with the following boundary conditions:
\begin{equation}\label{eq:boundary_cond}
\begin{cases}
	\textbf{u}^-
	=
	0
	& 
	\text{on } z = 0 \\
	\textbf{u}^- \cdot \textbf{t}^-
	=
	\textbf{u}^+ \cdot \textbf{t}^-
	& 
	\text{on } z = h^-(t,x) \\
	\textbf{u}^- \cdot \textbf{n}^-
	=
	\textbf{u}^+ \cdot \textbf{n}^+
	& 
	\text{on } z = h^-(t,x) \\
	\textbf{n}^- (\textbf{T}^- - \textbf{T}^+) \cdot \textbf{n}^-
	=
	\sigma^- \kappa^-
	&
	\text{on } z = h^-(t,x) \\
	\textbf{n}^- (\textbf{T}^- - \textbf{T}^+) \cdot \textbf{t}^-
	=
	0 
	&
	\text{on } z = h^-(t,x) \\
	(\textbf{T}^+ \textbf{n}^+) \cdot \textbf{n}^+
	=
	\sigma^+ \kappa^+ 
	&
	\text{on } z = h^+(t,x) \\
	(\textbf{T}^+ \textbf{n}^+) \cdot \textbf{t}^+
	=
	0
	&
	\text{on } z = h^+(t,x) \\
	h_t^{\pm} + u^{\pm} h_x^{\pm}
	=
	v^{\pm}
	& 
	\text{on } z = h^{\pm}(t,x).
\end{cases}
\end{equation}
Here, $\sigma^\pm$ is the constant surface tension and $\kappa^\pm$ is the mean curvature of the fluid-fluid interface $z=h^-(t,x) > 0$ and the fluid-air interface $z=h^+(t,x) > 0$, respectively.
%surface and $h^\pm = h^\pm(t,x) > 0$ the film height of the lower fluid ($-$) and the upper fluid ($+$), respectively. 
The first boundary condition is the so-called no-slip condition, saying that at the solid boundary the adhesive forces are stronger than the cohesive forces. Conditions $\eqref{eq:boundary_cond}_2$--$\eqref{eq:boundary_cond}_3$ guarantee that the tangential velocity as well as the normal velocity of the two immiscible fluids coincide at the fluid-fluid interface $z = h^-(t,x)$. The normal stress balance condition $\eqref{eq:boundary_cond}_4$ ensures that the jump in normal component of the stress across the fluid-fluid interface balances the curvature force per unit length. We have a similar condition at the fluid-air interface $z=h^+(t,x)$, taking into account that there is no hydrodynamic force exerted at the interface from above by the air, cf. $\eqref{eq:boundary_cond}_6$.
Since we assume to have constant surface tension $\sigma^\pm$ at both interfaces $z=h^\pm(t,x)$, there do not appear any gradients of $\sigma^\pm$, whence we assume in $\eqref{eq:boundary_cond}_5$ and $\eqref{eq:boundary_cond}_7$, respectively, that there is no jump in the tangential components of the stress at the interfaces. 
Finally, $\eqref{eq:boundary_cond}_8$ is the so-called kinematic boundary condition, saying that particles which are on the surface stay at the surface when time evolves.

\medskip
\subsection{Non-dimensional Navier--Stokes system}
In a next step we convert the above introduced Navier--Stokes system to non-dimensional form. To this end, we denote by $L$ the characteristic length and by $H$ the characteristic height of the fluid film. In addition, we introduce the parameter $\eps = \frac{H}{L} > 0$ and we assume that the film height is very small compared to its length, i.e. $H \ll L$ or equivalently $\eps \ll 1$. 
More precisely, we introduce dimensionless variables and unknowns
\begin{equation*}
	\begin{cases}
	\bar{x}=\frac{x}{L}, \quad \bar{z}=\frac{z}{H}, \quad \bar{t}=\eps^3 \frac{t}{\tau_0}, 
	\quad \bar{h}^\pm = \frac{h^\pm}{H} & \\
	\bar{u}^\pm = u^\pm  \frac{\tau _0}{L \eps^3}, \quad
	\bar{v}^\pm = v^\pm \frac{\tau_0}{L \eps^4}, \quad 
	%m^- \bar{\mu}^- = \mu^-, \quad
	\bar{\mu}^\pm\bigl(\tau \left\|\bar{\textbf{D}}^\pm\right\|\bigr)
	= \frac{1}{\mu_0^\pm}\mu^\pm\bigl(\tau_{\text{char}} \left\|\textbf{D}^\pm\right\|\bigr), & \\
	\tau = \frac{\tau_{\text{char}}}{\tau_0}, \quad
	\bar{p}^\pm = \tau_0 \frac{p^\pm}{\mu_0^\pm \eps}, & \\
	{m = \frac{\mu_0^+}{\mu_0^-}}, \quad s^\pm = \sigma^\pm \frac{\tau_0}{L \mu_0^\pm}, \quad \text{Re} = \frac{\rho^\pm \eps^3 L^2}{\mu_0^\pm \tau_0}.
	\end{cases}
\end{equation*}
Here $\tau_0$ is the macroscopic time scale of the system, $\tau_{\text{char}}$ is the characteristic time scale of the non-Newtonian fluids and $\mu_0^\pm$ is the characteristic viscosity of the respective fluid. Note that, due to the smallness of the parameter $\eps$, the scaling of the velocity field indicates that the horizontal velocity is relatively large compared to the vertical velocity. Moreover, the scaling of the pressure indicates that viscous forces are dominant, i.e. we are in the regime of so-called creeping flows.
Finally, $s^\pm$ is the non-dimensional surface tension and $\text{Re}$ denotes the Reynolds number. Note that the non-dimensional quantities $\tau, m$ and $\text{Re}$ do not have a subscript $\pm$.

In these dimensionless variables the Navier--Stokes system reads
\begin{equation*}
\begin{cases}
    \eps^2 \text{Re} \left( \bar{u}^\pm_{\bar{t}} 
    + \bar{u}^\pm \ \bar{u}^\pm_{\bar{x}}
    + \bar{v}^\pm\ \bar{u}^\pm_{\bar{z}} \right)
    = -\bar{p}^\pm_{\bar{x}} 
    + \left(\bar{\mu}^\pm \left( \eps^2\bar{u}^\pm_{\bar{x}\bar{x}} +\bar{u}^\pm_{\bar{z}\bar{z}}\right)
    +2\eps^2 \bar{\mu}^\pm_{\bar{x}}\bar{u}^\pm_{\bar{x}}
    + \bar{\mu}^\pm_{\bar{z}}\left( \bar{u}^\pm_{\bar{z}}
    + \eps^3 \bar{v}^\pm_{\bar{x}}\right) \right) 
    & 
    \text{in } \bar{\Omega}_{\pm}(\bar{t})
    \\
    \eps^5 \text{Re} \left( \bar{v}^\pm_{\bar{t}} 
    +  \bar{u}^\pm \ \bar{v}^\pm_{\bar{x}}
    + \bar{v}^\pm\ \bar{v}^i_{\bar{z}} \right)
    = - \bar{p}^\pm_{\bar{z}} 
    + \eps^2 \left(\bar{\mu}^\pm \left(\eps^2 \bar{v}^\pm_{\bar{x}\bar{x}} +\bar{v}^\pm_{\bar{z}\bar{z}}\right)
    +2  \bar{\mu}^\pm_{\bar{z}}\bar{v}^\pm_{\bar{z}}
    +  \bar{\mu}^\pm_{\bar{x}}  \left( \bar{u}^\pm_{\bar{z}} + \eps^2 \bar{v}^\pm_{\bar{x}} \right)\right) 
    & 
    \text{in } \bar{\Omega}_{\pm}(\bar{t})
    \\
    \bar{u}^\pm_{\bar{x}} +\bar{v}^\pm_{\bar{z}} 
    = 0
    &
    \text{in } \bar{\Omega}_{\pm}(\bar{t}),
    \end{cases}
\end{equation*}
where the rescaled domains $\bar{\Omega}_{\pm}(\bar{t})$ are given by
\begin{equation*}
\begin{cases}
	\bar{\Omega}_{-}(\bar{t}) = \left\{(\bar{x},\bar{z}) \in \R^2;\ \bar{x}\in \bar{\Omega},\ 0 < \bar{z} < \bar{h}^-(\bar{t},\bar{x})\right\},&
	\\
	\bar{\Omega}_{+}(\bar{t}) = \left\{(\bar{x},\bar{z}) \in \R^2;\ \bar{x}\in \bar{\Omega},\ \bar{h}^-(\bar{t},\bar{x}) < \bar{z} < \bar{h}^+(\bar{t},\bar{x}) \right\}.&
\end{cases}
\end{equation*}
%Taking the limit $\eps$ to zero, we end up with the equations
%\begin{equation*}
%\begin{cases}
%    -\bar{p}^\pm_{\bar{x}}
%    +(\bar{\mu}^\pm_{\bar{z}}\bar{u}^\pm_{\bar{z}})_{\bar{z}}
%    =0
%    & 
%    \text{in } \bar{\Omega}_{\pm}(\bar{t})  
%    \\
%    \bar{p}^\pm_{\bar{z}}
%    =0
%    & 
%    \text{in } \bar{\Omega}_{\pm}(\bar{t})
%    \\
%    \bar{u}^\pm_{\bar{x}} +\bar{v}^\pm_{\bar{z}} 
%    = 0
%    & 
%	\text{in } \bar{\Omega}_{\pm}(\bar{t}),
%\end{cases}
%\end{equation*}
Moreover, the no-slip condition and the kinematic boundary condition in dimensionless variables are
\begin{equation*}
\begin{cases}
    \bar{\textbf{u}}^- 
    =0 
    & 
    \text{on } \bar{z} = 0 
    \\
    \bar{\textbf{u}}^- 
    =\bar{\textbf{u}}^+ 
    & 
    \text{on } \bar{z} = \bar{h}^-(\bar{t},\bar{x}) 
    \\
    \bar{h}^\pm_{\bar{t}} + \bar{u}^\pm \bar{h}^\pm_{\bar{x}} 
    = \bar{v}^\pm 
    &
    \text{on } \bar{z}=\bar{h}^\pm(\bar{t},\bar{x}).
\end{cases}
\end{equation*}
Finally, with the relations
\begin{equation*}
	\textbf{n}^\pm = \frac{1}{\sqrt{1+|h^\pm_x|^2}} 
	\begin{pmatrix}
	-h^\pm_x \\
	1
	\end{pmatrix},
	\quad
	\textbf{t}^\pm = \frac{1}{\sqrt{1+|h^\pm_x|^2}} 
	\begin{pmatrix}
	1 \\
	h^\pm_x
	\end{pmatrix}
	\quad 
	\text{and}
	\quad
	\kappa^\pm
	=
	\frac{h^\pm_{xx}}{\left(1+|h^\pm_x|^2\right)^{\scriptstyle 3/2}}
\end{equation*}
for the outer pointing unit normal, the unit tangent and the mean curvature of the interfaces, 
we may write the rescaled normal-stress boundary condition at the fluid-fluid interface $\bar{z}=\bar{h}^-(\bar{t},\bar{x})$ and the fluid-air interface $\bar{z}=\bar{h}^+(\bar{t},\bar{x})$, respectively, as\footnote{For the sake of brevity we use the notation $A^\pm ]^+_- = A^+ - A^-$ in the following lines.}
\begin{equation*}
\begin{cases}
	m^\pm \bar{p}^\pm\bigr]^+_-
	\left( 1+ \eps^2|\bar{h}^-_{\bar{x}}|^2 \right) 
	-2 \eps^2 
	m^\pm \bar{\mu}^\pm 
	\left(\eps^2 |\bar{h}^-_{\bar{x}}|^2 \bar{u}^\pm_{\bar{x}}
	- \bar{h}^-_{\bar{x}} 
	\left(\bar{u}^\pm_{\bar{z}}+ \eps^2 \bar{v}^\pm_{\bar{x}} \right) 
	+\eps^2 \bar{v}^\pm_{\bar{z}}\right)\bigr]^+_-
	=
	\frac{\sigma^- \bar{h}^-_{\bar{x}\bar{x}}}{L \tau_0\sqrt{1+\eps^2 |\bar{h}^-_{\bar{x}}|^2}} &
	\\
	-\left( 1+ \eps^2|\bar{h}^+_{\bar{x}}|^2 \right) \bar{p}^+
	+ 2 \eps^4 \bar{\mu}^+ |\bar{h}^+_{\bar{x}}|^2 \bar{u}^+_{\bar{x}}
	-2 \eps^2 \bar{\mu}^+ \bar{h}^+_{\bar{x}} 
	\left( \bar{u}^+_{\bar{z}} + \eps^2 \bar{v}^+_{\bar{x}} \right)
	+ 2 \eps^2 \bar{\mu}^+ \bar{v}^+_{\bar{z}} 
	= \frac{\sigma^+ \bar{h}^+_{\bar{x}\bar{x}}}{Lm^+ \tau_0 \sqrt{1+\eps^2 |\bar{h}^+_{\bar{x}}|^2}}, &
\end{cases}
\end{equation*}
and the rescaled tangential-stress balance condition at the respective interfaces $\bar{z}=\bar{h}^\pm(\bar{t},\bar{x})$ as 
\begin{equation*}
\begin{cases}
	-2 \tau_0 \eps^4 \bar{h}^-_{\bar{x}} \left(
	m^\pm \bar{\mu}^\pm \left( \bar{v}^\pm_{\bar{z}} - \bar{u}^\pm_{\bar{x}} \right)\right)\bigr]^+_-
	- \tau_0 \eps^2 \left( 1- \eps^2 |\bar{h}^-_{\bar{x}}|^2 \right) 
	\left(
	m^\pm \bar{\mu}^\pm \left(\bar{u}^\pm_{\bar{z}} + \eps^2 \bar{v}^\pm_{\bar{x}} \right)\right)\bigr]^+_-= 0
	\\
	2 \eps^4 \tau_0 m^+ \bar{\mu}^+ \left( \bar{v}^+_{\bar{z}} -\bar{u}^+_{\bar{x}} \right)
	+ \left( 1-\eps^2|\bar{h}^+_{\bar{x}}|^2 \right) \tau_0 m^+ \eps^2 \bar{\mu}^+
	\left( \bar{u}^+_{\bar{z}} + \eps^2 \bar{v}^+_{\bar{x}} \right) = 0. &
\end{cases}
\end{equation*}
For convenience of the reader from now on we omit the bars in our notation. 
%%------------------------------
%%------------------------------
%%------------------------------
\medskip

\subsection{Thin-film approximation}
We now exploit the assumption that both liquid films are of positive but very small height. This means that we consider the formal asymptotic limit $\eps \searrow 0$. For a rigorous justification of the so-called lubrication approximation for a single Newtonian thin film, the Stokes and Hele--Shaw flow the reader is referred to \cite{GO, GP, MP}.
In our setting of a relatively slow flow, we obtain the reduced Navier--Stokes equations
\begin{equation}\label{eq:Navier-Stokes_limit}
\begin{cases}
    -p^\pm_x + \left( \mu^\pm(\tau |u_z^\pm|) u^\pm_z \right)_z 
    = 0 
    &
    \text{in } \Omega^\pm(t) %\label{eq:pres_x} 
    \\
    p^\pm_z
    =0 
    & 
    \text{in } \Omega^\pm(t) %\label{eq:pres_z} 
    \\
    u^\pm_x+v^\pm_z
    =0 
    & 
    \text{in } \Omega^\pm(t) %\label{eq:incomp} 
\end{cases}
\end{equation}
and the boundary conditions become
\begin{equation} \label{eq:boundary_cond_limit}
\begin{cases}
    \textbf{u}^-
    =0 
    & 
    \text{on } z=0 %\label{bc:bt} 
    \\
    \textbf{u}^-
    =\textbf{u}^+ 
    & 
    \text{on } z=h^-(t,x) %\label{bc:interface} 
    \\
    h^\pm_t+u^\pm h^\pm_x 
    = v^\pm 
    & 
    \text{on } z=h^\pm(t,x) %\label{bc:kin} 
    \\
    m p^+ - p^- 
    = s^- h^-_{xx} 
    & \text{on } z=h^-(t,x) %\label{stress_bal1} 
    \\
    p^+ 
    = - s^+ h^+_{xx} 
    & 
    \text{on } z=h^+(t,x) %\label{stress_bal2} 
    \\
    \mu^-(\tau |u_z^-|) u^-_z 
    = 
    m \mu^+(\tau |u_z^+|) u^+_z 
    & 
    \text{on } z=h^-(t,x) %\label{stress_bal3} 
    \\
    u^+_z
    = 0 
    & \text{on } z=h^+(t,x).  %\label{stress_bal4}
\end{cases}
\end{equation}
From this system we can explicitly determine the pressure $p^\pm$ and the velocity $\textbf{u}^\pm$ of both fluids. Then we are left with a system of evolution equations for the film heights 
\begin{equation*}
	f(t,x)=h^-(t,x)
	\quad
	\text{and}
	\quad 
	g(t,x)=h^+(t,x)-h^-(t,x), \quad t > 0,\, x \in \Omega.
\end{equation*}
Indeed, we integrate $\eqref{eq:Navier-Stokes_limit}_2$ and use the relations $\eqref{eq:boundary_cond_limit}_6$ and $\eqref{eq:boundary_cond_limit}_7$ to obtain the equations
\begin{equation*}
	\begin{cases}
    p^+(t,x) = - s^+ (f+g)_{xx} &\\
    p^-(t,x) = -m s^+(f+g)_{xx} - s^- f_{xx} &
    \end{cases}
\end{equation*}
for the pressure. Moreover, 
integrating $\eqref{eq:Navier-Stokes_limit}_1$ with respect to $z$ and using again the boundary conditions $\eqref{eq:boundary_cond_limit}_6$ and $\eqref{eq:boundary_cond_limit}_7$ yields
\begin{equation}\label{eq:mu_u_z}
\begin{cases}
    \mu^+(\tau |u_z^+|) u^+_z = s^+ (f+g)_{xxx} (f+g-z) & \\
    \mu^-(\tau |u_z^-|) u^-_z = m s^+ (f+g)_{xxx} (f+g-z) 
    + s^- f_{xxx} (f-z). &
\end{cases}
\end{equation}
In order to derive evolution equations for the two interfaces we assume at this point that the shear stresses are increasing functions of the shear rates, i.e. the functions $s \mapsto \mu^\pm(|s|) s,\ s \in \R$, are increasing. This allows us to resolve the equations \eqref{eq:mu_u_z} for $u_z^\pm$. Indeed, under this assumption we can define functions $\psi^\pm$ such that $\psi^\pm\bigl(\mu^\pm(|s|)s\bigr) = s,\ s \in \R$. Dropping for the moment the $t$-dependence in the notation, this implies that
\begin{equation*}%\label{eq:u_z}
\begin{cases}
	u^+_z(x,z) = \frac{1}{\tau}\psi^+\bigl(\tau s^+ (f+g)_{xxx} (f+g-z)\bigr) & \\
	u^-_z(x,z) = \frac{1}{\tau}\psi^-\bigl(\tau m s^+ (f+g)_{xxx} (f+g-z) 
	+ \tau s^- f_{xxx} (f-z)\bigr) &
\end{cases}
\end{equation*}
and integration with respect to $z$ yields
\begin{equation}\label{eq:u}
\begin{cases}
	u^+(x,z) = \frac{1}{\tau} \displaystyle \int_{f(x)}^z \psi^+\bigl(\tau s^+ (f(x)+g(x))_{xxx} (f(x)+g(x)-\xi)\bigr)\, d\xi & \\
	u^-(x,z) = \frac{1}{\tau} \displaystyle\int_0^z\psi^-\bigl(\tau m s^+ (f(x)+g(x))_{xxx} (f(x)+g(x)-\xi) 
	+ \tau s^- f_{xxx}(x) (f(x)-\xi)\bigr)\, d\xi &
\end{cases}
\end{equation}
for all $t > 0$ and $x \in \Omega$.
Thanks to the conservation-of-mass equation one can therewith determine $v^\pm$. The only unknowns which are yet to be determined are the film heights $f$ and $g$. For this purpose we observe that, using the conservation of mass equation $\eqref{eq:Navier-Stokes_limit}_3$, the kinematic boundary conditions $\eqref{eq:boundary_cond_limit}_3$ may be rewritten as
\begin{equation*}
	\begin{cases}
	f_t + \left(\displaystyle \int_0^{f(x)} u^-(x,z)\, dz\right)_x = 0, & t > 0,\ x \in \Omega 
	\\
	g_t + \left(\displaystyle \int_{f(x)}^{f(x)+g(x)} u^+(x,z)\, dz\right)_x = 0, & t > 0,\ x \in \Omega.
	\end{cases}
\end{equation*}
Inserting the expressions \eqref{eq:u} for the horizontal velocity yields
\begin{equation*}
	\begin{cases}
		f_t + \frac{1}{\tau} \left(\displaystyle \int_0^{f(x)} \displaystyle\int_0^z \psi^-\bigl(\tau m s^+ (f(x)+g(x))_{xxx} (f(x)+g(x)-\xi) 
		+ \tau s^- f_{xxx}(x) (f(x)-\xi)\bigr)\, d\xi dz\right)_x = 0 & 
		\\
		g_t + \frac{1}{\tau} \left(\displaystyle \int_{f(x)}^{f(x)+g(x)} \displaystyle \int_{f(x)}^z  \psi^+\bigl(\tau s^+ (f(x)+g(x))_{xxx} (f(x)+g(x)-\xi)\bigr)\, d\xi\, dz\right)_x = 0 &
	\end{cases}	
\end{equation*}
for all $t > 0$ and $x \in \Omega$. Using Fubini's theorem, we find that
\begin{equation*}
	\begin{split}
		\int_0^{f(x)} \int_0^z & \psi^-\bigl(\tau m s^+ (f(x)+g(x))_{xxx} (f(x)+g(x)-\xi) 
		+ \tau s^- f_{xxx}(x) (f(x)-\xi)\bigr)\, d\xi dz
		\\
		&= \int_0^{f(x)} \int_\xi^{f(x)} \psi^-\bigl(\tau m s^+ (f(x)+g(x))_{xxx} (f(x)+g(x)-\xi) 
		+ \tau s^- f_{xxx}(x) (f(x)-\xi)\bigr)\, dz\, d\xi  
		\\
		&= \int_0^{f(x)} (f(x)-\xi)\ \psi^-\bigl(\tau m s^+ (f(x)+g(x))_{xxx} (f(x)+g(x)-\xi) 
		+ \tau s^- f_{xxx}(x) (f(x)-\xi)\bigr)\, d\xi 
		\\
		&= \int_0^{f(x)} r\ \psi^-\bigl(\tau m s^+ (f(x)+g(x))_{xxx} (g(x)+r) 
		+ \tau s^- f_{xxx}(x) r\bigr)\, dr
		\\
		&= f(x)^2 \int_0^1 y \psi^-\bigl(\tau m s^+ (f(x)+g(x))_{xxx} (g(x)+y f(x)) 
		+ \tau s^- f_{xxx}(x) y f(x)\bigr)\, dy
	\end{split}
\end{equation*}
and similarly
\begin{equation*}
	\begin{split}
		\int_{f(x)}^{f(x)+g(x)} & \int_{f(x)}^z \psi^+\bigl(\tau s^+ (f(x)+g(x))_{xxx} (f(x)+g(x)-\xi)\bigr)\, d\xi\, dz
		\\
		&= \int_{f(x)}^{f(x)+g(x)} \int_\xi^{f(x)+g(x)} \psi^+\bigl(\tau s^+ (f(x)+g(x))_{xxx} (f(x)+g(x)-\xi)\bigr)\, dz\, d\xi
		\\
		&= \int_{f(x)}^{f(x)+g(x)} (f(x)+g(x)-\xi)\, \psi^+\bigl(\tau s^+ (f(x)+g(x))_{xxx} (f(x)+g(x)-\xi)\bigr)\, d\xi
		\\
		&=
		\int_{f(x)}^{f(x)+g(x)} r\, \psi^+\bigl(\tau s^+ (f(x)+g(x))_{xxx}\, r\bigr)\, dr
		\\
		&=
		\int_0^1 \bigl(yg(x)^2+f(x)g(x)\bigr)\, \psi^+\bigl(\tau s^+ (f(x)+g(x))_{xxx}\, (yg(x)+f(x))\bigr)\,dy.
	\end{split}
\end{equation*}
Consequently, the evolution equations for $f$ and $g$ may further be rewritten as
\begin{equation*}
\begin{cases}
	f_t + \frac{1}{\tau} \left(\displaystyle f(x)^2 \int_0^1 y \psi^-\bigl(\tau m s^+ (f(x)+g(x))_{xxx} (g(x)+y f(x)) 
	+ \tau s^- f_{xxx}(x) y f(x)\bigr)\, dy \right)_x = 0 & 
	\\
	g_t + \frac{1}{\tau} \left(\displaystyle \int_0^1 \bigl(yg(x)^2+f(x)g(x)\bigr)\, \psi^+\bigl(\tau s^+ (f(x)+g(x))_{xxx}\, (yg(x)+f(x))\bigr)\,dy\right)_x = 0 &
\end{cases}	
\end{equation*}
for all $t > 0,\, x \in \Omega$. Changing the time variable via
\begin{equation*}
	\hat{t} = \frac{1}{\tau} t
\end{equation*}
and dropping the $\hat{\phantom{u}}$ in the notation for convenience, we end up with the system 
\begin{equation*}
\begin{cases}
	f_t + \left(\displaystyle f(x)^2 \int_0^1 y \psi^-\bigl(\tau m s^+ (f(x)+g(x))_{xxx} (g(x)+y f(x)) 
	+ \tau s^- f_{xxx}(x) y f(x)\bigr)\, dy \right)_x = 0 & 
	\\
	g_t + \left(\displaystyle \int_0^1 \bigl(yg(x)^2+f(x)g(x)\bigr)\, \psi^+\bigl(\tau s^+ (f(x)+g(x))_{xxx}\, (yg(x)+f(x))\bigr)\,dy\right)_x = 0 &
\end{cases}	
\end{equation*}
of evolution equations for the film heights $f = f(t,x)$ and $g=g(t,x)$, where $t > 0$ and $x \in \Omega$.
%%----------------------------
%%----------------------------
%%----------------------------
\medskip

\subsection{An Ellis fluid on top of a Newtonian fluid}
In this subsection we consider the case in which the lower film $\Omega_-(t)$ is filled by a Newtonian fluid, i.e. $\mu^- \equiv 1$, while the upper film $\Omega_+(t)$ is filled by an Ellis fluid. The dynamic viscosity of the upper fluid is thus implicitly given by the constitutive law
\begin{equation}\label{eq:Ellis}
	\frac{1}{\mu^+(\tau|u^+_z|)}=
	\frac{1}{\mu^+_0} \left(1 + \left|\frac{\mu^+(\tau|u^+_z|)\, u^+_z}{\tau^+_{1/2}}\right|^{p-1}\right) \quad \text{with } p \geq 1.
\end{equation}
Here $\mu^+(\tau|u^+_z|)\, u^+_z$ is the shear stress,
$\mu_0^+$ is the viscosity at zero shear, $\tau_{1/2}^+$ is the shear-stress at which the viscosity drops down to $\mu_0^+/2$, and $p \geq 1$ is the so-called flow-behaviour exponent.
For $p = 1$ the fluid is Newtonian, for $p > 1$ it is shear-thinning. We now derive the evolution equations for the film heights $f$ and $g$ in this special regime. First, using in \eqref{eq:mu_u_z} that $\mu^-=1$, we have
\begin{equation*}
	u_z^-(t,x,z) = m s^+ (f+g)_{xxx} (f+g-z) + s^- f_{xxx}(f-z), \quad t > 0,\, x \in \Omega,
\end{equation*}
for all $z \in (0,f)$.
In this case the function $\psi^-$ is the identity. The horizontal velocity of the lower fluid is thus given by
\begin{equation*}
    u^-(t,x,z) = m s^+ (f+g)_{xxx} \left(fz+gz-\frac{z^2}{2} \right) 
    + s^- f_{xxx} \left(fz-\frac{z^2}{2}\right).
\end{equation*}
Using the Ellis constitutive law \eqref{eq:Ellis} for the upper fluid in \eqref{eq:mu_u_z}, and introducing the notation
\begin{equation*}
	\Phi(d) = |d|^{p-1} d, \quad d\in \R,
\end{equation*}
we find that
\begin{align*}
    u^+_z(t,x,z) &= 
    \frac{s^+}{\mu^+_0} 
    \left(1 + \left|\frac{\mu^+(\tau|u^+_z|)\, u^+_z}{\tau^+_{1/2}}\right|^{p-1}\right)
    (f+g)_{xxx} (f+g-z) \\
%    &= \frac{s^2}{\mu^2_0} (f+g)_{xxx} (f+g-z)
%    + \frac{s^2}{\mu^2_0 |\tau^2_{1/2}|^{\beta-1}}
%    \left|\mu^2(|u^2_z|)u^2_z\right|^{\beta-1}
%    (f+g)_{xxx} (f+g-z) \\
    &= \frac{s^+}{\mu^+_0} (f+g)_{xxx} (f+g-z)
    + \frac{|s^+|^p}{\mu^+_0 |\tau^+_{1/2}|^{p-1}}
    \Phi\left((f+g)_{xxx} (f+g-z)\right)
\end{align*}
for all $t > 0,\, x\in \Omega$ and $z \in (f,f+g)$.
Integrating this equation from $f$ to $z$ and using the interface condition $\eqref{eq:boundary_cond_limit}_2$ implies that the horizontal velocity of the upper fluid is given by
\begin{align*}
    u^+(t,x,z) &= m s^+ (f+g)_{xxx} \left( \frac{f^2}{2} + f g \right) + \frac{s^-}{2} f_{xxx} f^2
    %\\&
    + \frac{s^+}{\mu^+_0} (f+g)_{xxx} \left(f z + g z -\frac{z^2}{2} -\frac{f^2}{2} - fg \right)
    \\ &\quad
    - \frac{|s^+|^p}{(p+1) \mu^+_0 |\tau^+_{1/2}|^{p-1}}
    \Phi\left((f+g)_{xxx}\right) \left(|f+g-z|^{p+1} - |g|^{p+1}\right)
\end{align*}
for all $t > 0, x \in \Omega$ and $z \in (f,f+g)$. Consequently, in the case in which an Ellis shear-thinning fluid is placed on a Newtonian fluid, we obtain that the evolution equations for $f$ and $g$ are given by
\begin{equation} \label{eq:system_0}
\hspace{-1cm}	\begin{cases}
		f_t + \left(m s^+ \left( \frac{f^3}{3}+ \frac{f^2 g}{2} \right) (f+g)_{xxx}
		+\frac{s^-}{3} f^3 f_{xxx}\right)_x = 0 & 
		\\
		g_t + 
		\left(\left(\frac{ms^+}{2} f^2g + ms^+ f g^2 + \frac{s^+}{3 \mu_0^+} g^3\right)(f+g)_{xxx} + \frac{s^-}{2} f^2 g f_{xxx} + C_p |g|^{p + 2} \Phi\bigl((f+g)_{xxx}\bigr)\right)_x = 0, & 
	\end{cases}
\end{equation}
for $t > 0,\, x \in \Omega$, where the constant $C_p$ is defined as
\begin{equation} \label{eq:C_p}
C_p = \frac{|s^+|^{p}}{(p+2) \mu_0^+ |\tau_{1/2}^+|^{p-1}}.
\end{equation}
We shall study this system on a finite interval $\Omega\subset \R$, equipped with positive initial conditions
\[
f_0(x),g_0(x)>0, \quad x \in \Omega,
\]
and the Neumann-type boundary conditions
\begin{equation}\label{eq:boundary} 
	f_x=f_{xxx}=0\quad \mbox{and}\quad g_x=g_{xxx}=0\qquad \mbox{on}\quad \partial \Omega.
\end{equation}             
Let us briefly comment on the main properties of the system \eqref{eq:system_0} of evolution equations. The equations for both film heights $f$ and $g$ are quasilinear equations of fourth order, stated in divergence form. Moreover, they are degenerate in $f$ and $g$, respectively. In the equation for the height $f$ of the Newtonian fluid film the dependence of the coefficients on the solution is smooth. Note that this is not the case in the equation for the height $g$ of the non-Newtonian fluid film. Recalling the definition of the function $\Phi$, we observe that, for flow behaviour exponents $p \in (1,2)$, the coefficients of the highest-order term depend only $(p-1)$-H\"older continuously on the third-order spatial derivatives $(f+g)_{xxx}$. This H\"older continuous dependence may be seen explicitly in Section \ref{sec:existence}, where we give up the divergence form in order to prove existence of strong solutions for short times.

%=============================================================================
%=============================================================================
%=============================================================================
\bigskip

\section{Local existence of strong solutions}\label{sec:existence}

In this section we prove existence of strong solutions to the problem
%\begin{equation}\tag{P}\label{eq:system}
%\begin{cases}
%	f_t + \left(m s^+ \left( \frac{f^3}{3}+ \frac{f^2 g}{2} \right) (f+g)_{xxx}
%	+\frac{s^-}{3} f^3 f_{xxx}\right)_x = 0, \quad t > 0,\ x \in \Omega & 
%	\\
%	g_t + 
%	\left(\left(\frac{ms^+}{2} f^2g + ms^+ f g^2 + \frac{s^+}{3 \mu_0^+} g^3\right)(f+g)_{xxx} + \frac{s^-}{2} f^2 g f_{xxx} + C_p |g|^{p + 2} \Phi\bigl((f+g)_{xxx}\bigr)\right)_x = 0,\ t > 0,\ x \in \Omega & 
%	\\
%	f_x=f_{xxx}=g_x=g_{xxx}=0, \quad t > 0,\ x \in \partial\Omega &
%	\\
%	f(0,x)=f_0(x),\ g(0,x)=g_0(x), \quad x \in \Omega &
%\end{cases}
%\end{equation}
%for short times, where the positive constant $C_p > 0$ is defined in \eqref{eq:C_p}.

\begin{flalign}\tag{P}\label{eq:system}
	\begin{cases}
		f_t + \left(m s^+ \left( \frac{f^3}{3}+ \frac{f^2 g}{2} \right) (f+g)_{xxx}
		+\frac{s^-}{3} f^3 f_{xxx}\right)_x = 0 
		\quad \text{in } \Omega_T & 
		\\
		g_t + 
		\left(\mkern-6mu\left(\mkern-3mu\frac{ms^+}{2} f^2g + ms^+ f g^2 + \frac{s^+}{3 \mu_0^+} g^3 \mkern-3mu\right)(f+g)_{xxx} + \frac{s^-}{2} f^2 g f_{xxx} + C_p |g|^{p + 2} \Phi\bigl((f+g)_{xxx}\bigr)\mkern-6mu\right)_{\mkern-6mu x}
		\mkern-3mu = 0\quad \text{in } \Omega_T & 
		\\
		f_x=f_{xxx}=g_x=g_{xxx}=0 \quad \text{on } \Gamma_T &
		\\
		f(0,x)=f_0(x),\ g(0,x)=g_0(x) \quad \text{for } x\in\Omega, &
	\end{cases}&&
\end{flalign}
with $\Omega_T=(0,T)\times\Omega$ and $\Gamma_T=(0,T)\times\partial\Omega$ for some possibly small $T>0$. The constant $C_p > 0$ is as defined in \eqref{eq:C_p}.
We set $u=(f,g)$ and recast \eqref{eq:system} as an abstract quasilinear Cauchy problem of the form
\begin{equation*}
	\begin{cases}
		u_t + \Acal(u)u = \Fcal(u), & t > 0 \\
		u(0) = u_0
	\end{cases}
\end{equation*}	
on a Banach space involving the boundary conditions.
Such a problem can be solved by means of abstract semigroup theory. We start by specifying the functional analytic setting we work in. The underlying space on which we study the system \eqref{eq:system} of evolution equations is the Hilbert space $L_2(\Omega)$. The operator $\Acal$ is a differential operator of fourth order which is to be defined such that the fourth-order derivative is in again in $L_2(\Omega)$ and such that the Neumann-type boundary conditions are incorporated in its domain. For this purpose, we introduce for $s \geq 0$ the $L_2$-based Bessel potential spaces $H^s(\Omega)$.
%Let us introduce suitable Banach spaces on which the system of evolution equations \eqref{eq:system} is going to be studied.
%For $s \geq 0$ we denote by $H^s(\Omega)$ the $L_2$-based Bessel potential spaces.
These are defined as restrictions of functions in $H^s(\R)$ to $\Omega$ with the norm
\[
%\|v\|_{H^s(\Omega)}:=\inf\{\|w\|_{H^s(\R)};\ w\in H^s(\Omega), w|_\Omega=v\},
\|v\|_{H^s(\Omega)}:=\inf\{\|w\|_{H^s(\R)};\ w\in H^s(\R), w|_\Omega=v\},
\]
where
\[
H^s(\R):=\{w\in L_2(\R);\ (1+|\cdot|^2)^\frac{s}{2}\hat w \in L_2(\R)\}.
\]
Here, $\hat w$ denotes the Fourier transformation of $w$.
Recall that whenever $s$ is a natural number, the space $H^s(\Omega)$ coincides with the classical Sobolev space
\[
W^s_2(\Omega):=\{v\in L_2(\Omega);\ \partial_x^n v\in L_2(\Omega)\mbox{ for all } n\leq s \}.
\]
Furthermore, 
\begin{equation*}
	H^s(\Omega) = [H^{s_0}(\Omega),H^{s_1}(\Omega)]_{\rho}, 
\end{equation*}
where $\rho\in (0,1)$, $0\leq s_0 < s_1 < \infty$, and $s = (1 - \rho) s_0 + \rho s_1$. That is, $H^s(\Omega)$ appears as the complex interpolation space
between $H^{s_1}(\Omega)$ and $H^{s_0}(\Omega)$. 
Moreover, in order to take the first-order and third-order Neumann-type boundary conditions into account, we introduce the Banach spaces
\begin{equation*}
	H^{4\rho}_B(\Omega)
	=
	\begin{cases}
	\left\{v \in H^{4\rho}(\Omega);\ v_x = v_{xxx} = 0 \text{ on }  \partial\Omega\right\}, & \frac{7}{2} < 4\rho \leq 4
	\\[1ex]
	\left\{v \in H^{4\rho}(\Omega);\ v_x = 0 \text{ on } \partial\Omega\right\}, & \frac{3}{2} < 4\rho \leq \frac{7}{2}
	\\[1ex] 
	H^{4\rho}(\Omega), & 0 \leq 4\rho \leq \frac{3}{2}.
	\end{cases}
\end{equation*}
%The space $H^4_B(\Omega)\times H^4_B(\Omega)$ plays a natural role in our analysis. 
For $4\rho \in (0,4)\setminus\{3/2,7/2\}$, the spaces $H^{4\rho}_B(\Omega)$ are closed linear subspaces of $H^{4\rho}(\Omega)$ and satisfy the interpolation property
\begin{equation} \label{eq:interpolation}
	H^{4\rho}_B(\Omega)
	=
	\bigl[L_2(\Omega),H^4_B(\Omega)\bigr]_{\rho}.
\end{equation}
%======================
Now, we are in a position to define the abstract quasilinear Cauchy problem corresponding to \eqref{eq:system} in this setting.
%This section is devoted to prove classical local existence for \eqref{eq:system} in $H_B^{4\theta}(\Omega)\times H_B^{4\theta}(\Omega)$ for $\theta>\frac{7}{8}$. Let us recast the system of evolution equations \eqref{eq:system} as an abstract Cauchy problem. For this purpose, we introduce the following notation. For $r > s > 1/2$ we define
%\begin{equation*}
%	\nu = \frac{3+s}{4} \quad \text{and} \quad \theta = \frac{3+r}{4}
%\end{equation*}
%such that $\theta > \nu > 7/8$. 
For this purpose, we introduce
for $u=(f,g)\in H^{4\rho}_B(\Omega)\times H^{4\rho}_B(\Omega)$ and $\rho>7/8$ the differential matrix operator
%For this purpose, we set $u=(f,g)$ and introduce
%the fourth-order (linear) matrix operator
\begin{equation}\label{eq:differential_op}
	\begin{cases}
		\mathcal{A}: H^{4\rho}_B(\Omega)\times H^{4\rho}_B(\Omega)\longrightarrow \mathcal{L}\bigl(H^4_B(\Omega)\times H^4_B(\Omega);L_2(\Omega)\times L_2(\Omega)\bigr) & \\
		\mathcal{A}(u)w:=A(u,u_{xxx})\partial_x^4 w, \quad 
		w \in H^4_B(\Omega)\times H^4_B(\Omega) &
%		u \in H^{4\rho}_B(\Omega)\times H^{4\rho}_B(\Omega) &
	\end{cases}
\end{equation}
of fourth order,
%\[
%\mathcal{A}: H^{4\rho}_B(\Omega)\times H^{4\rho}_B(\Omega)\to \mathcal{L}\bigl(H^4_B(\Omega)\times H^4_B(\Omega),L_2(\Omega)\times L_2(\Omega)\bigr)
%\]
%by 
%\begin{equation*}
%	\mathcal{A}(u,u_{xxx})w:=A(u,u_{xxx})\partial_x^4 w, \quad u \in H^{4\rho}_B(\Omega)\times H^{4\rho}_B(\Omega),
%\end{equation*}
where the matrix $A(u,u_{xxx})$ is defined by
\begin{equation*}
	A(u,u_{xxx})=
	\begin{pmatrix}
	a_{11}(u) & a_{12}(u) \\
	a_{21}(u,u_{xxx}) & a_{22}(u,u_{xxx})
	\end{pmatrix}
\end{equation*}
with components
\begin{equation*}
	\begin{cases}
		a_{11}(u)=\frac{ms^+ + s^-}{3}f^3+\frac{ms^+}{2}f^2g, &
		\\
		a_{12}(u)=m s^+ \left(\frac{f^3}{3}+ \frac{f^2 g}{2} \right), &
		\\
		a_{21}(u,u_{xxx})=\frac{m s^+ + s^-}{2}f^2 g + ms^+ f g^2 + \frac{s^+}{3\mu_0^+}g^3 + C_p g^{p+2} \Phi'\left((f+g)_{xxx}\right), &
		\\
		a_{22}(u,u_{xxx})=\frac{m s^+}{2} f^2 g + m s^+ f g^2 + \frac{s^+}{3\mu_0^+} g^3 + C_p g^{p+2} \Phi'\left((f+g)_{xxx}\right). &
	\end{cases}
\end{equation*}
%The domain of the differential operator $\Acal$ is given by
%\begin{equation*}
%	\mathcal{D}(\Acal) = \left\{w \in L_2(\Omega)\times L_2(\Omega);\ \Acal(u)w \in L_2(\Omega)\times L_2(\Omega)\ \forall u \in H^{4\rho}_B(\Omega)\times H^{4\rho}_B(\Omega),\ u_x=u_{xxx}=0 \text{ on } \partial\Omega\right\}
%\end{equation*}
%{\color{orange} GB: We should recheck the domain. The condition $u_x=u_{xxx}=0$ is already included in $H_B^{4\rho}$. Is it a typo and you wanted to write: $w_x=w_{xxx}=0$? Do we actually have to state the domain of $\Acal$? Before, we say $\Acal$ is an operator form $H_B^{4\rho}$ to $\mathcal{L}(H_B^{4\rho},L_2)$. Maybe we can simply leave it like this?}
%\begin{align*}
%	a_{11}(u)&=\frac{ms^+ + s^-}{3}f^3+\frac{ms^+}{2}f^2h, \\
%	a_{12}(u)&=m s^+ \left(\frac{f^3}{3}+ \frac{f^2 h}{2} \right), \\
%	a_{21}(u,u_{xxx})&=\frac{m s^+ + s^-}{2}f^2 h + ms^+ f h^2 + \frac{s^+}{3\mu_0^+}h^3 + C h^{\beta+2} \Phi'\left((f+h)_{xxx}\right), \\
%	a_{22}(u,u_{xxx})&=\frac{m s^+}{2} f^2 h + m s^+ f h^2 + \frac{s^+}{3\mu_0^2+} h^3 + C h^{\beta+2} \Phi'\left((f+h)_{xxx}\right).
%\end{align*}
The nonlinear lower-order terms of \eqref{eq:system} are comprised in the map 
\begin{equation}\label{eq:rhs}
	\begin{cases}
		\mathcal{F}: H^{4\rho}_B(\Omega)\times H^{4\rho}_B(\Omega)\longrightarrow L_2(\Omega)\times L_2(\Omega) & \\
		\Fcal(u) = 
		\begin{pmatrix}
			F_1(u,u_x,u_{xxx})\\
			F_2(u,u_x,u_{xxx})
		\end{pmatrix},
		\quad u \in H^{4\rho}_B(\Omega)\times H^{4\rho}_B(\Omega),
	\end{cases}
\end{equation}
%$\mathcal{F}: H^{4\rho}_B(\Omega)\times H^{4\rho}_B(\Omega)\to L_2(\Omega)\times L_2(\Omega)$ defined by
%\begin{equation*}
%	\Fcal(u,u_x,u_{xx},u_{xxx}) = 
%	\begin{pmatrix}
%		\Fcal_1(u,u_x,u_{xx},u_{xxx})\\
%		\Fcal_2(u,u_x,u_{xx},u_{xxx})
%	\end{pmatrix},
%	\quad u \in H^{4\rho}_B(\Omega)\times H^{4\rho}_B(\Omega),
%\end{equation*}
where 
\begin{equation*}
	\begin{cases}
	F_1(u,u_x,u_{xxx}) = m s^+
	\left( f^2 f_x + f g f_x +\frac{1}{2} f^2 g_x \right)
	(f+g)_{xxx} + s^- f^2 f_x f_{xxx}, & \\
	F_2(u,u_x,u_{xxx}) = m s^+ \left( f g f_x 	+\frac{1}{2} f^2 g_x 
	+ g^2 f_x + 2 f g g_x 
	+ \frac{1}{m \mu_0^+}g^2 g_x \right)
	(f+g)_{xxx} \\
	\phantom{F_2(u,u_x,u_{xx},u_{xxx}) = } +s^-\left( f g f_x + \frac{1}{2} f^2 g_x \right) 
	f_{xxx}
	+ C (p+2) g^{p+1} 
	\Phi\left( (f+g)_{xxx} \right) g_x. &
	\end{cases}
\end{equation*}
%\begin{align*}
%	F_1(u,u_x,u_{xx},u_{xxx})&= m s^+
%	\left( f^2 f_x + f g f_x +\frac{1}{2} f^2 g_x \right)
%	(f+g)_{xxx} + s^- f^2 f_x f_{xxx},\\
%	F_2(u,u_x,u_{xx},u_{xxx})&= m s^+ \left( f g f_x +\frac{1}{2} f^2 g_x 
%	+ g^2 f_x + 2 f g g_x 
%	+ \frac{1}{m \mu_0^+}g^2 g_x \right)
%	(f+g)_{xxx}\\
%	&\quad +s^-\left( f g f_x + \frac{1}{2} f^2 g_x \right) 
%	f_{xxx}
%	+ C (p+2) g^{p+1} 
%	\Phi\left( (f+g)_{xxx} \right) g_x.
%\end{align*}
With this notation, we can write \eqref{eq:system} as an abstract quasilinear Cauchy problem 
\begin{equation} \tag{CP}\label{eq:Cauchy}
\begin{cases}
	u_t + \Acal(u)\, u = \Fcal(u), & t > 0
	\\
	u(0) = u_0
\end{cases}
\end{equation}
on $L_2(\Omega)\times L_2(\Omega)$ with initial datum $u_0=(f_0,g_0)$. 
We call \eqref{eq:Cauchy} parabolic, if $-\Acal(u)$ is the infinitesimal generator of an analytic semigroup on $L_2(\Omega)\times L_2(\Omega)$ with domain $H^{4}_B(\Omega)\times H^{4}_B(\Omega)$. We denote the set of negative infinitesimal generators of analytic semigroups on a Banach space $E_0$ by  $\Hcal(E_1;E_0)$, where $E_1\subset E_0$ is the domain of the generator. The space $\Hcal(E_1;E_0)$ is equipped with the norm $\left\|\cdot\right\|_{\Lcal(E_1;E_0)}$.

A classical result of Amann \cite[Theorem 12.1]{A:1993} yields the existence of a unique strong solution of \eqref{eq:Cauchy} for positive initial data in appropriate function spaces if the maps 
\begin{equation*}
	\begin{cases}
	\mathcal{A}: H^{4\rho}_B(\Omega)\times H^{4\rho}_B(\Omega)\longrightarrow \mathcal{L}\bigl(H^4_B(\Omega)\times H^4_B(\Omega);  L_2(\Omega)\times L_2(\Omega)\bigr) &\\
	\mathcal{F}: H^{4\rho}_B(\Omega)\times H^{4\rho}_B(\Omega)\longrightarrow L_2(\Omega)\times L_2(\Omega) &
	\end{cases}
\end{equation*} 
are Lipschitz continuous and $\Acal(u)\in \mathcal{H}\bigl(H^4_B(\Omega)\times H^4_B(\Omega);  L_2(\Omega)\times L_2(\Omega)\bigr)$ for every $u\in H^{4\rho}_B(\Omega)\times H^{4\rho}_B(\Omega)$. This is the case when $p=1$, i.e. when the upper fluid is Newtonian, or when $p\geq 2$. If $p\in (1,2)$, then $\Acal$ and the derivative of $\Fcal$ are merely $(p-1)$-H\"older continuous and problem \eqref{eq:system} does not fall into the frame of \cite[Theorem 12.1]{A:1993}. However, on the a priori cost of uniqueness of solutions, we may apply a recently developed result on existence of classical solutions of quasilinear Cauchy problems of the form \eqref{eq:Cauchy}, where the dependence of the functionals $\Acal$ and $\Fcal$ on its coefficients is assumed to be only H\"older continuous \cite[Theorem 4.2]{LM:2020}: 

\begin{theorem}\label{thm:abstract}
Suppose that $(E_0,E_1)$ is a densely and compactly injected Banach couple. 
For $\eta \in (0,1)$, we denote by $E_\eta:=[E_0,E_1]_\eta$ the real or the complex interpolation functor. 
Let $0<\beta<\alpha\leq 1$ and $\sigma\in (0,\alpha-\beta)$. Let $\mu\in (0,1)$ and assume that 
\[
\mathcal{A} \in C^\mu_{\text{loc}}\bigl(E_\beta; \mathcal{H}(E_1,E_0)\bigr)\qquad \mbox{and}\qquad \mathcal{F}\in C^\mu_{\text{loc}}(E_\beta;E_0).
\]
Then, for each $u_0\in E_\alpha$ there exists $T>0$ such that the quasilinear Cauchy problem
\begin{equation*}
\begin{cases}
	u_t + \Acal(u)\, u = \Fcal(u), & t > 0,
	\\
	u(0) = u_0
\end{cases}
\end{equation*}
possesses a solution
\[
u\in C^{\mu\sigma}\bigl((0,T];E_1\bigr)\cap C^{1+\mu\sigma}\bigl((0,T]; E_0\bigr)\cap C\bigl([0,T];E_\alpha\bigr).
\]
Moreover, for any $\alpha^\prime\in (\beta+\sigma,\alpha)$ the solution satisfies
\[
u\in C^\sigma\left([0,T]; E_{\alpha^\prime-\sigma}\right)\cap C^{\alpha^\prime-\beta}\left([0,T];E_\beta\right).
\]
\end{theorem}

Equipped with this theorem for quasilinear parabolic Cauchy problems with H\"older continuous dependence, we prove the existence of local strong solutions of \eqref{eq:Cauchy}. 
%======================

\begin{theorem}[Local strong solutions of \eqref{eq:system}] \label{thm:existence}
Let $\theta > \rho > 7/8$. Given $u_0=(f_0,g_0) \in  H^{4\theta}_B(\Omega)\times H^{4\theta}_B(\Omega)$
with $f_0,g_0>0$, for each $p > 1$ there exists a positive time $T > 0$ and a solution $u=(f,g)$ of \eqref{eq:Cauchy} on $[0,T]$ with initial datum $u_0$, satisfying the regularity
	\begin{equation*}
	\begin{split}
		(f,g) \in\ & C\bigl([0,T]; H^{4\theta}_B(\Omega)\times H^{4\theta}_B(\Omega)\bigr)
		\cap
		C\bigl((0,T];H^4_B(\Omega)\times H^4_B(\Omega)\bigr)
		\cap
		C^1\bigl((0,T];L_2(\Omega)\times L_2(\Omega)\bigr).
	\end{split}
	\end{equation*}
	{In addition, the solution is H\"older continuous in the sense that
	\begin{equation*}
		u \in C^{\nu}\bigl([0,T];H^{4\rho}_B(\Omega)\times H^{4\rho}_B(\Omega)\bigr)
	\end{equation*}
	for $\nu \in (0,\theta-\rho)$ and}
	and it satisfies
	\begin{equation*}
		f(t,x), g(t,x) > 0, \quad (t,x) \in [0,T]\times\bar{\Omega}.
	\end{equation*}
\end{theorem}

%======================

%To see that Theorem \ref{thm:existence} holds true, observe that \eqref{eq:system} can be classified into the abstract setting of \cite[Theorem 4.2]{LM:2020}.
% 
%====================== 

\begin{proof}
In order to apply Theorem \ref{thm:abstract}, we identify 
\[
E_0=L_2(\Omega)\times L_2(\Omega)\qquad \mbox{and}\qquad E_1=H^4_B(\Omega)\times H^4_B(\Omega).
\]
Moreover, for $\theta > \rho > 7/8$, we have that
\begin{equation*}
	E_\theta(\Omega) = H^{4\theta}_B(\Omega)\times H^{4\theta}_B(\Omega) \quad \text{and} \quad E_\rho(\Omega) = H^{4\rho}_B(\Omega)\times H^{4\rho}_B(\Omega).
\end{equation*}
It is well known that the embeddings
\[
	H^4_B(\Omega)\times H^4_B(\Omega) \longhookrightarrow
	H^{4\theta}_B(\Omega)\times H^{4\theta}_B(\Omega) \longhookrightarrow 
	H^{4\rho}_B(\Omega)\times H^{4\rho}_B(\Omega) \longhookrightarrow
	L_2(\Omega)\times L_2(\Omega)
\]
are dense and compact, c.f. for instance \cite[Chapter I.2.5]{A:1995}. 
%Moreover, $E_\theta=H^{4\theta}(\Omega)\times H^{4\theta}(\Omega)$. 
Due to the choices of $\theta>\rho > 7/8$, we have in addition that $H^{4\theta}_B(\Omega)\times H^{4\theta}_B(\Omega)\hookrightarrow H^{4\rho}_B(\Omega)\times H^{4\rho}_B(\Omega) \hookrightarrow C^3(\bar{\Omega})\times C^3(\bar{\Omega})$, and the maps
\begin{equation*}
	\begin{cases}
		\mathcal{A}: H^{4\rho}_B(\Omega)\times H^{4\rho}_B(\Omega)\longrightarrow \mathcal{L}\bigl(H^4_B(\Omega)\times H^4_B(\Omega);L_2(\Omega)\times L_2(\Omega)\bigr) & \\
		\mathcal{F}: H^{4\rho}_B(\Omega)\times H^{4\rho}_B(\Omega)\longrightarrow L_2(\Omega)\times L_2(\Omega) &
	\end{cases}
\end{equation*}
are locally $(p-1)$-H\"older continuous. Due to the degeneracy of \eqref{eq:system}, we certainly lose parabolicity % of $\Acal(u)$ 
{of \eqref{eq:Cauchy}}
when one of the two film heights tends to zero. The idea is that starting with a strictly positive initial datum, the operator $-\Acal(u(t))$ is the infinitesimal generator of an analytic semigroup as long as $u(t)\in H^{4\rho}_B(\Omega)\times H^{4\rho}_B(\Omega)$ is strictly positive. To make the proof rigorous, we fix $\eps>0$ and introduce an extended operator $\Acal_\eps$ defined by $\Acal_\eps(u)=A_\eps(u,u_{xxx})\partial_x^4$, where
\[
A_\eps(u,u_{xxx})=\begin{pmatrix} a_{11}(u_\eps) & a_{12}(u_\eps)\\ a_{21}(u_\eps,u_{xxx}) & a_{22}(u_\eps,u_{xxx})\end{pmatrix},\qquad u_\eps=(f_\eps,g_{\eps}):=\bigl(\max(f,\eps),\max(g,\eps)\bigr).
\]
Note that the map $u\mapsto u_\eps$ is locally Lipschitz continuous and we do not extend the coefficients with respect to $u_{xxx}$. Hence, we still have that the map
\[
\Acal_\eps:H^{4\rho}_B(\Omega)\times H^{4\rho}_B(\Omega) \longrightarrow \mathcal{L}\bigl(H^4_B(\Omega)\times H^4_B(\Omega);L_2(\Omega)\times L_2(\Omega)\bigr)
\]
is locally $(p-1)$-H\"older continuous. We prove that for any $u\in H^{4\rho}_B(\Omega)\times H^{4\rho}_B(\Omega)$, the operator $-\Acal_\eps(u)$ is the infinitesimal generator of an analytic semigroup on $L_2(\Omega)\times L_2(\Omega)$, which implies, in view of Theorem \ref{thm:existence}, that there exists a time $T_\eps>0$ such that the extended problem 
\begin{equation} \label{eq:Cauchy_ex}
\begin{cases}
	u_t + \Acal_\eps(u)\, u = \Fcal(u), & t > 0
	\\
	u(0) = u_0
\end{cases}
\end{equation}
has a solution
\begin{equation*}
	u(\eps;\,\cdot) \in C\bigl([0,T_\eps]; H^{4\theta}_B(\Omega)\times H^{4\theta}_B(\Omega)\bigr)
	\cap
	C^\nu\bigl([0,T_\eps]; H^{4\rho}_B(\Omega)\times H^{4\rho}_B(\Omega)\bigr)
	\cap 
	C^1\bigl((0,T_\eps];L_2(\Omega)\times L_2(\Omega)\bigr),
\end{equation*}
where $\nu \in (0,\theta-\rho)$.
%	\begin{equation*}
%	\begin{split}
%		\textcolor{magenta}{u} \in\ & C\bigl([0,T_\eps]; H^{4\theta}_B(\Omega)\times H^{4\theta}_B(\Omega)\bigr)
%		\cap
%		C\bigl((0,T_\eps];H^4_B(\Omega)\times H^4_B(\Omega)\bigr)
%		\cap
%		C^1\bigl((0,T_\eps];L_2(\Omega)\times L_2(\Omega)\bigr).
%	\end{split}
%	\end{equation*}
Eventually, we show that $\Acal_\eps(u(\eps;t))=\Acal(u(\eps;t))$ for $t\in [0,T_\eps]$, which implies the existence of a strictly positive solution to the original problem \eqref{eq:Cauchy} on the time interval $[0,T_\eps]$.\medskip

%\subsection*{Generator of an analytic semigroup}
\noindent\textsc{Generator of an analytic semigroup. }
Let $u\in H^{4\rho}_B(\Omega)\times H^{4\rho}_B(\Omega)$ and fix $\eps>0$. Note that then 
\[A_\eps(u,u_{xxx})\in C\bigl(\bar{\Omega},\R^{2\times2}\bigr)
\] since the entries satisfy $a_{1i}(u_\eps)\in C^3(\bar{\Omega},\R)$ and $a_{2i}(u_\eps, u_{xxx})\in C(\bar{\Omega},\R)$, $i=1,2$, due to the embedding $H^{4\rho}_B(\Omega)\subset C^3(\bar{\Omega})$ for $\rho>\frac{7}{8}$. In order to show that $-\mathcal{A}_\eps(u)$ is the infinitesimal generator of an analytic semigroup on $L_2(\Omega)\times L_2(\Omega)$ with domain $H^4_B(\Omega)\times H^4_B(\Omega)$, we verify that 
\begin{itemize}
\item[(a)] \textit{$\Acal_\eps(u)$ is normally elliptic}. That is,
\begin{equation*}
	\sigma(a^{\pi}_\eps(x,\xi))\subset \{z\in \C\mid \mbox{Re}\, z >0\}\qquad \mbox{for all}\qquad (x,\xi) \in \bar{\Omega}\times\{-1,1\},
\end{equation*}
where $a^{\pi}_\eps$ denotes the principal symbol of $\Acal_\eps(u)$.
\item[(b)] \textit{$\Acal_\eps(u)$ satisfies the Lopatinskii--Shapiro condition}. That is, zero is the only exponentially decaying solution of the ordinary differential equation
\begin{equation}\label{eq:ODE_LS}
	\begin{cases}
		\left(\lambda
		+a^{\pi}_\eps(x,\ii\partial_t)\right)\,\phi(t)=0, 
		&  t>0 \\
		b^\pi(\ii \partial_t)\,\phi(0)=0 &
	\end{cases}
\end{equation}
for $\lambda\in\{z\in\C \mid \text{Re}\ z>0\}$,
where $b^\pi$ denotes the principal boundary symbol.
\end{itemize}
Then, it follows from \cite[Theorem 4.1 \& Remark 4.2(b)]{A:1993} that $-\Acal_\eps(u)$ generates an analytic semigroup on $L_2(\Omega)\times L_2(\Omega)$. 
{In order to check conditions (a) and (b) note first that the} principal symbol of $\Acal_\eps(u)$ can be written as
\begin{equation*}
	a^\pi_\eps (x,\xi)=A_\eps(u,u_{xxx}) \xi^4,
	\quad (x,\xi)\in \bar{\Omega}\times \R.
\end{equation*}
Hence, $\sigma(a^\pi_\eps(x,\xi)) \subset \{ z\in\C :\text{Re}\, z>0 \}$ if and only if the eigenvalues of the matrix $A_\eps(u,u_{xxx})$ have strictly positive real part. Recall that
\[
	A_\eps(u,u_{xxx})=\begin{pmatrix} a_{11}(u_\eps) & a_{12}(u_\eps) \\ a_{21}(u_\eps, u_{xxx}) &  a_{22}(u_\eps, u_{xxx})\end{pmatrix},
\]
where $a_{ij}>0$ for $1\leq i,j\leq 2$ and the entries $a_{11}$ and $a_{21}$ can be written as
\begin{equation*}
	\begin{cases}
		a_{11}(u_\eps)=a_{12}(u_{\eps})+ \frac{s^-}{3}f_\eps^3 &\\
		a_{21}(u_\eps, u_{xxx})=a_{22}(u_\eps,u_{xxx})+\frac{s^-}{2}f_\eps^2g_\eps.&
	\end{cases}
\end{equation*}
%Notice that then 
{Using these representations it can be seen that $A_\eps(u,u_{xxx})$ has a strictly positive determinant. Indeed, it holds that}
\begin{align*}
\det(A_\eps(u,u_{xxx}))&=a_{11}(u_\eps)a_{22}(u_\eps,u_{xxx})-a_{12}(u_\eps)a_{21}(u_\eps,u_{xxx})\\
&=\left(a_{12}(u_{\eps})+ \frac{s^-}{3}f_\eps^3\right)a_{22}(u_\eps,u_{xxx})-a_{12}(u_\eps)\left(a_{22}(u_\eps,u_{xxx})+\frac{s^-}{2}f_\eps^2g_\eps\right)\\
&=\frac{s^-}{3}f_\eps^3 a_{22}(u_\eps,u_{xxx})-\frac{s^-}{2}f_\eps^2g_\eps a_{12}(u_\eps)\\
&=\frac{m s^- s^+}{12} f_\eps^4 g_\eps^2
	+\frac{s^- s^+}{9 \mu_0^+} f_\eps^3 g_\eps^3
	+ \frac{C_p s^-}{3}f_\eps^3 g_\eps^{p+2}
	\Phi'\left((f+g)_{xxx}\right)\\
	&>0.
\end{align*}
Since the characteristic polynomial of $A_\eps(u,u_{xxx})$ is given by
\begin{equation*}
	\text{det}(\lambda-A_\eps(u,u_{xxx}))
	= \lambda^2 - \left(a_{11}(u_\eps)+a_{22}(u_\eps,u_{xxx})\right) \lambda
	+a_{11}(u_\eps)a_{22}(u_\eps,u_{xxx}) - a_{12}(u_\eps) a_{21}(u_\eps,u_{xxx}),
\end{equation*}
the eigenvalues of $A_\eps(u,u_{xxx})$ are determined by
\begin{align}\label{eq:eigenvalues_A}
	\lambda_{\pm}&= \frac{a_{11}(u_\eps)+a_{22}(u_\eps,u_{xxx})}{2} 
	\pm \sqrt{\frac{\bigl(a_{11}(u_\eps)+a_{22}(u_\eps,u_{xxx})\bigr)^2}{4}
	-\det(A_\eps(u,u_{xxx}))}.
%&	= \frac{a_{11}(u_\eps)+a_{22}(u_\eps),u_{xxx}}{2} 
%	\pm \sqrt{\frac{(a_{11}(u_\eps)-a_{22}(u_\eps),u_{xxx})^2}{4}
%	+a_{12}(u_\eps)a_{21}(u_\eps,u_{xxx})}.
\end{align}
Using that $\det(A_\eps(u,u_{xxx}))>0$, we deduce that the eigenvalues of $A_{\eps}(u,u_{xxx})$ have strictly positive real part, which proves that condition (a) is fulfilled.

Next, we verify the Lopatinskii--Shapiro condition (b), where in our case the principal boundary symbol $b^\pi$ is given by
\begin{equation*}
	b^\pi (\xi)=\ii
	\begin{pmatrix}
	\xi & 0 \\
	0 & \xi \\
	-\xi^3 & 0 \\
	0 & -\xi^3
\end{pmatrix}.
\end{equation*}
Observing that $a^\pi_\eps(x,i\partial_t)=A_\eps(u,u_{xxx})\partial_t^4$, problem \eqref{eq:ODE_LS} reads
\begin{equation*}
	\begin{cases}
		\lambda\phi(t)
		+A_\eps(u,u_{xxx})\partial_t^4\phi(t)=0, 
		\quad  t>0 & \\
		\partial_t\phi(0)=\partial_t^3\phi(0)=0.&
	\end{cases}
\end{equation*}
Since $\text{det}(A_\eps(u,u_{xxx}))>0$, we know that $A_\eps(u,u_{xxx})$ is invertible and we have
\begin{equation*}
	\partial_t^4\phi(t)=-\lambda A_\eps(u,u_{xxx})^{-1}\,\phi(t).
\end{equation*}
Moreover, we can compute that
\[
\frac{\bigl(a_{11}(u_\eps)+a_{22}(u_\eps,u_{xxx})\bigr)^2}{4}
	-\det(A_\eps(u,u_{xxx}))=\frac{\bigl(a_{11}(u_\eps)-a_{22}(u_\eps,u_{xxx})\bigr)^2}{4}
	+a_{12}(u_\eps)a_{21}(u_\eps,u_{xxx})>0,
\]
which implies by \eqref{eq:eigenvalues_A} that the eigenvalues $\lambda_{\pm}$ of $A_\eps(u,u_{xxx})$ are real and $\lambda_-<\lambda_+$.
%
%\textcolor{red}{Ich dachte, das ist vllt etwas viel, weil die Bedingung ja auf der vorherigen Seite steht... Man muss halt das boundary symbol noch hinschreiben...Vorschlag:}
%Concerning the Lopatinskii--Shapiro condition (b), we consider the ordinary differential equation
%\begin{equation}\label{eq:ODE_LS}
%	\left(\lambda
%	+a^{\pi}_\eps(x,\ii\partial_t)\right)\,\phi(t)=0, 
%	\quad  t>0, \qquad
%	b^\pi(\ii \partial_t)\,\phi(0)=0, 
%\end{equation}
%where $b_\pi$ denotes the principal boundary symbol given by
%\begin{equation*}
%	b^\pi (\xi)=\ii
%	\begin{pmatrix}
%	\xi & 0 \\
%	0 & \xi \\
%	-\xi^3 & 0 \\
%	0 & -\xi^3
%\end{pmatrix}.
%\end{equation*}
%We aim to show that $\varphi \equiv 0$ is the only exponentially decaying solution of \eqref{eq:ODE_LS} for $\lambda \in \{z\in \C \mid \mbox{Re}\, z>0\}$.
Consequently, there exists an invertible matrix $U \in \C^{2\times 2}$  such that we can diagonalise  $A_\eps(u,u_{xxx})$ via
\begin{equation*}
	\text{diag}(\lambda_- ,\lambda_+) = U^{-1} A_\eps(u,u_{xxx}) \ U.
\end{equation*}
Setting $\psi(t)=U\phi(t)$, we get
\begin{equation*}
	\partial_t^4 \psi(t) = U \partial_t^4\phi(t) = -\lambda\, \bigl(U^{-1}A_\eps(u,u_{xxx})\ U\bigr)^{-1}\psi(t)=-\lambda \bigl(	\text{diag}(\lambda_- ,\lambda_+)\bigr)^{-1}\psi(t).
\end{equation*}
Hence, the system decouples into two independent equations
\begin{equation}\label{eq:ODE}
	\partial_t^4 \psi_k(t) = -\frac{\lambda}{\lambda_k}  \psi_k(t), \quad k=1, 2,
\end{equation}
where $\lambda_1=\lambda_-$ and $\lambda_2=\lambda_+$. The general solution to the ordinary differential equation \eqref{eq:ODE} of fourth order is given by
\[
\psi_k(t)=\sum_{n=1}^4c_ne^{\Lambda_{k,n} t},\quad t>0,
\]
where $\{\Lambda_{k,n};\, n=1,\ldots,4\}$ solve $\Lambda_k^4=-\frac{\lambda}{\lambda_k}$. Since $\lambda_k>0$ and $\mbox{Re} \lambda>0$, we have that $-\frac{\lambda}{\lambda_k}$ does not belong to the imaginary axis. {Hence, its fourth roots are given by two complex numbers with strictly positive real parts (say $\Lambda_{k,1}$ and $\Lambda_{k,2}$) and two complex numbers with strictly negative real parts (say $\Lambda_{k,3}$ and $\Lambda_{k,4}$).} Now, if $\psi_k$ is a an exponentially decaying solution, then
\[
\psi_k(t)=c_3e^{\Lambda_{k,3}t}+c_4e^{\Lambda_{k,4}t}, \qquad k=1,2.
\]
The initial conditions for $\phi$ imply that $\partial_t\psi(0)=\partial_t^3\psi(0)=0$ and hence $\psi=0$. We conclude by invertability of $U$ that  $\phi\equiv 0$ is the only exponentially decaying solution of \eqref{eq:ODE_LS} for $\lambda \in \C$ 
%with $\mbox{Re}\, \lambda>0$, which proves the Lopatinskii--Shapiro condition (b). 
with $\mbox{Re}\, \lambda>0$. This proves the Lopatinskii--Shapiro condition (b).
 
In view of  \cite[Theorem 4.1 \& Remark 4.2(b)]{A:1993} we have shown that %$-\Acal_\eps(u,u_{xxx})$ 
$-\Acal_\eps(u)$
generates an analytic semigroup on $L_2(\Omega)\times L_2(\Omega)$. Hence, due to Theorem \ref{thm:abstract}, there exists a time $T_\eps>0$ such that the extended problem \eqref{eq:Cauchy_ex} admits a solution\footnote{{Clearly the solution has even better regularity properties for strictly positive times, cf. Theorem \ref{thm:abstract}. However, here we mention only the regularity properties that we need in the following.}} 
\begin{equation*}
	u(\epsilon;\cdot) \in C\bigl([0,T_\eps]; H^{4\theta}_B(\Omega)\times H^{4\theta}_B(\Omega)\bigr)
	\cap
	C^\nu\bigl([0,T_\eps]; H^{4\rho}_B(\Omega)\times H^{4\rho}_B(\Omega)\bigr)
	\cap 
	C^1\bigl((0,T_\eps];L_2(\Omega)\times L_2(\Omega)\bigr),
\end{equation*}
where $\nu \in (0,\theta-\rho)$. 
\medskip
%\[u(\eps) \in C\left([0,T_\eps], H_B^{4\theta}(\Omega)\times H_B^{4\theta}(\Omega)\right)\cap C\left((0,T_\eps];H^4_B(\Omega)\times H^4_B(\Omega)\right)\cap C^1\left((0,T_\eps];L_2(\Omega)\times L_2(\Omega)\right).
%\] \medskip
%\subsection*{Positivity and solution of the original problem}

\noindent\textsc{Positivity and solution of the original problem. }
In order to prove the existence of a local strong solution to the original problem \eqref{eq:Cauchy}, we verify that for sufficiently small times $T$, the solution $u(\eps;\cdot)$ of the extended problem is strictly positive and $\Acal_\eps(u(\eps;t))=\Acal(u(\eps;t))$ for $t\in[0,T]$. 
Indeed, if $\min_{x\in \bar{\Omega}}(f_0,g_0)>2\eps$, then there exists a constant $c>0$ such that
\[
\min_{t\in [0,T_\eps]}\min_{x\in \bar{\Omega}}\bigl(f(\eps;t,x),g(\eps;t,x)\bigr)\geq 2\eps-cT_\eps^{\nu},
\]
where we used the notation $u(\eps;t,x)=(f(\eps;t,x),g(\eps;t,x))$. {That is, as long as $T_\eps<T^\ast:=\left(\frac{\eps}{c}\right)^{\frac{1}{\nu}}$, both components $f(\eps;\cdot)$ and $g(\eps;\cdot)$ of $u(\eps;\cdot)$ are strictly larger than $\eps$.}
%That is,  $u(\eps)$ is component wise strictly larger than $\eps$ for $T^*<\left(\frac{\eps}{c}\right)^{\frac{1}{\theta}}$.} 
By definition of $\Acal_\eps$, we infer that
\[
\Acal(u(\eps;t))=\Acal_\eps(u(\eps;t))
\qquad \text{for }t\in[0,T^*).
\]
Therefore, $u(\eps;\cdot)$ is a solution to the original problem \eqref{eq:Cauchy} on $[0,T^*)$ and the proof is complete.
\end{proof}

%=============================================================================
\bigskip

\begin{remark}
Similar as in \cite[Section 5]{LM:2020}, Theorem \ref{thm:existence} can be proved in other functional analytic settings. Examples are the fractional Sobolev spaces or the little H\"older spaces.
\end{remark}

%=============================================================================
%=============================================================================
%=============================================================================
\bigskip

\section{Uniqueness of strong solutions}\label{sec:uniqueness}
This section is concerned with uniqueness of solutions to \eqref{eq:system} for all flow-behaviour exponents $p>1$. Recall that for $p \geq 2$ the  the map $u\mapsto\Acal(u)$, introduced in \eqref{eq:differential_op}, as well as the right-hand side $\Fcal$, defined in \eqref{eq:rhs}, are at least Lipschitz continuous. In this case, the traditional results \cite[Theorem 12.1]{A:1993} and \cite[Theorem III.4.6.3]{E:1969} of \textsc{Amann} and \textsc{Eidel'man}, respectively,  are applicable and yield existence and uniqueness of solutions to \eqref{eq:system} by a contraction argument in the underlying fixed-point problem. For flow-behaviour exponents $p \in (1,2)$ problem \eqref{eq:system} lacks this Lipschitz property, i.e. the  map $u\mapsto\Acal(u)$ is merely $(p-1)$-H\"older continuous. Existence may, however, be proved by compactness, cf. Theorem \ref{thm:existence}.
%=============================================================================
%=============================================================================
%=============================================================================
% Recall that the existence theorem of the previous section yields uniqueness of strong solutions only for flow behaviour exponents $p \geq 2$ in the constitutive law for the viscosity function $\mu^+$ describing the rheological behaviour of the upper fluid. This is due to a lack of regularity in the coefficient matrix of the differential operator.
In this section we prove for $p > 1$ that, if two positive solutions of \eqref{eq:system} coincide initially, then they coincide as long as they exist. This is done by defining a suitable energy functional and by exploiting the divergence form of the system of differential equations.

In the following  it is  convenient to rewrite \eqref{eq:system} as a system of evolution equations for $f$ and $f+g$ in divergence form, that is
\begin{equation} \label{eq:system_new}
\begin{cases}
	f_t + \bigl(  p_{11}(f,g)f_{xxx} + p_{12}(f,g)(f+g)_{xxx}\bigr)_x =0 & 
	\\
	(f+g)_t +\bigl( p_{21}(f,g)f_{xxx} + p_{22}(f,g)(f+g)_{xxx} + C_p |g|^{p +2}\Phi((f+g)_{xxx}) \bigr)_x =0 ,
\end{cases}
\end{equation}
where
\begin{equation*}
    \begin{cases}
        p_{11}(f,g) = \frac{s^-}{3}f^3 & \\
        p_{12}(f,g) = ms^+\left(\frac{1}{3}f^3+\frac{1}{2}f^2g\right) & \\
        p_{21}(f,g) = p_{11}(f,g)+ \frac{s^-}{2}f^2g &\\
        p_{22}(f,g) = p_{12}(f,g)+ms^+\left(f g^2 +\frac{1}{2}f^2g\right)+\frac{s^+}{3{\mu_0^+}}g^3. &
    \end{cases}
\end{equation*}
% \begin{align*}
% p_{11}(f,g)&= \frac{s^-}{3}f^3\\
% p_{12}(f,g)&= ms^+\left(\frac{1}{3}f^3+\frac{1}{2}f^2g\right)\\
% p_{21}(f,g)&=p_{11}(f,g)+ \frac{s^-}{2}f^2g\\
% p_{22}(f,g)&=p_{12}(f,g)+ms^+\left(fg^2 +\frac{1}{2}f^2g\right)+\frac{s^+}{3{\mu_0^+}}g^3
% \end{align*}
We define the functional
\begin{equation*}
	E(f,g)
	:=
	\frac{1}{2} \int_{\Omega} \left|(f+g)_x\right|^2 + \frac{s^-}{ms^+} \left|f_x\right|^2\, dx
\end{equation*}
and observe that $E=E(f,g)$ is an energy functional in the sense that it decreases along smooth solutions of \eqref{eq:system_new}. More precisely, for solutions $(f,g)$ of \eqref{eq:system_new} in the sense of Theorem \ref{thm:existence} the following statement holds true.

 %\textcolor{magenta}{CL: smooth sol. or sol. in the sense of the ex. thm?!}

%\textcolor{cyan}{CL: Somehow we mention the energy identity twice in the lemma. What if we write it as follows?!}

\begin{lemma}[Energy functional]\label{lem:Energy}
If $(f,g)$ is a solution of \eqref{eq:system_new} the sense of Theorem \ref{thm:existence} on $[0,T]$ with initial datum $(f_0,g_0)$, then it satisfies the energy identity
\begin{equation}\label{eq:energy_id}
   E(f,g)(t)+\int_0^t D(f,g)(\tau)\,d\tau = E(f_0,g_0),\quad t\in [0,T], 
\end{equation}
where the dissipative term $D(f,g)$ is defined by
\begin{align*}
    D(f,g)&=\int_{\Omega}\left(p_{21}(f,g)f_{xxx}+p_{22}(f,g)(f+g)_{xxx}+C_p |g|^{p +2}\Phi((f+g)_{xxx})\right)(f+g)_{xxx}\, dx\\
	&+\int_{\Omega} \frac{s^-}{ms^+}\left(p_{11}(f,g)f_{xxx}+p_{12}(f,g)(f+g)_{xxx}\right)f_{xxx}\, dx
\end{align*}
and satisfies $D(f,g)(t)\geq 0$ for all $t \in (0,T)$.
\end{lemma}
% \begin{lemma}[Energy functional]\label{lem:Energy}
% If $(f,g)$ is a solution of \eqref{eq:system_new} the sense of Theorem \ref{thm:existence} on $[0,T]$ with initial datum $(f_0,g_0)$, then
% \[
% \frac{d}{dt}E(f,g)(t)=-D(f,g)(t),\quad t \in (0,T),
% \]
% where
% \begin{align*}
% D(f,g)&:=\int_{\Omega}\left(p_{21}(f,g)f_{xxx}+p_{22}(f,g)(f+g)_{xxx}+C_p |g|^{p +2}\Phi((f+g)_{xxx})\right)(f+g)_{xxx}\, dx\\
% 	&+\int_{\Omega} \frac{s^-}{ms^+}\left(p_{11}(f,g)f_{xxx}+p_{12}(f,g)(f+g)_{xxx}\right)f_{xxx}\, dx.
% \end{align*}
% We have that $D(f,g)(t)\geq 0$ on $(0,T)$ and
% \[
% E(f,g)(t)+\int_0^tD(f,g)(\tau)\,d\tau = E(f_0,g_0),\quad t\in [0,T].
% \]
% \end{lemma}
\begin{proof}
%Notice fist that the regularity of a solution $(f,g)$ of \eqref{eq:system_new} due to Theorem \ref{thm:existence} guarantees that the map $t \mapsto E(f,g)(t)$ is differentiable on $(0,T)$. Indeed, $f_x,g_x \in C\left((0,T);H^1_0(\Omega)\times H^1(\Omega)\right)\cap C^1\left((0,T);H^{-1}(\Omega)\times H^{-1}(\Omega)\right) $, where $H^1_0$ denotes the space of functions belonging to $H^1(\Omega)$ with zero boundary condition, and $H^{-1}(\Omega)$ its dual space. Let $s,t\in (0,T)$, then
%\begin{align*}
%E(f,g)(t)-E(f,g)(s)
%&=2\int_s^t \langle\frac{d}{dt}(f_x+g_x)(\tau),(f_x+g_x)(\tau) \rangle\, d\tau +2\frac{\sigma^+}{\sigma_-}\int_s^t \langle\frac{d}{dt}f_x(\tau),f_x(\tau) \rangle\, d\tau,
%\end{align*}
%where $\langle \cdot, \cdot \rangle$ denotes the dual pairing between $H^1_0(\Omega)$ and $H^{-1}(\Omega)$.  Noting that $\langle\frac{d}{dt}f_x,f_x \rangle$ and $\langle\frac{d}{dt}g_x,g_x \rangle$ are continuous on $(0,T)$, we infer that the difference quotient $\frac{1}{t-s}\left(E(f,g)(t)-E(f,g)(s)\right)$ is finite.

Let $(f,g)$ be a solution of \eqref{eq:system_new} as in Theorem \ref{thm:existence}. Then we have in particular 
\begin{equation*}
	f_x,g_x \in 
	C\left((0,T);H^1_0(\Omega)\times H^1_0(\Omega)\right)
	\cap 
	C^1\left((0,T);H^{-1}(\Omega)\times H^{-1}(\Omega)\right), 
\end{equation*}
where $H^1_0$ denotes the space of functions belonging to $H^1(\Omega)$ with zero boundary condition, and $H^{-1}(\Omega)$ its dual space. It follows from \cite[Chapter 5.9, Theorem 3]{E:2010} that the map 
$t \mapsto E(f,g)(t)$ is absolutely continuous and we can compute its derivative as
\begin{align*}
	\frac{d}{dt} E(f,g)&=
	\int_\Omega (f+g)_{x} (f+g)_{xt} + \frac{s^-}{ms^+} f_x f_{xt}\, dx
	\\
	&=- \int_\Omega (f+g)_{xx} (f+g)_t + \frac{s^-}{ms^+} f_{xx} f_t\, dx
	\\
	&=-\int_{\Omega}\left(p_{21}(f,g)f_{xxx}+p_{22}(f,g)(f+g)_{xxx}+C_p |g|^{p +2}\Phi((f+g)_{xxx})\right)(f+g)_{xxx}\, dx\\
	&-\int_{\Omega} \frac{s^-}{ms^+}\left(p_{11}(f,g)f_{xxx}+p_{12}(f,g)(f+g)_{xxx}\right)f_{xxx}\, dx\\
	&=-D(f,g),
\end{align*}
where we use integration by parts and the boundary conditions $f_{xxx}=g_{xxx}=0$ on $\partial \Omega$.
%\textcolor{blue}{Alternatively, we could use the following:}\\
%For $\phi\in C\left((0,T);H^1_0(\Omega)\times H^1_0(\Omega)\right)\cap C^1\left((0,T);H^{-1}(\Omega)\times H^{-1}(\Omega)\right)$ it holds that
%\begin{equation*}
%	\frac{\norm{\phi(t)}_{L^2}^2-\norm{\phi(s)}_{L^2}^2}{t-s}
%	= \left\langle\frac{\phi(t)-\phi(s)}{t-s}, \ 
%	\phi(t)+\phi(s)\right\rangle 
%	\rightarrow 2 \langle \phi'(t),\phi(t) \rangle
%	\quad \text{as } s\rightarrow t.
%\end{equation*}
Using Lemma \ref{lem:A} with $A=f_{xxx}$ and $B=g_{xxx}$,
we complete squares such that
\begin{align}\label{eq:D}
\begin{split}
D(f,g)&=\int_\Omega C_p g^{p+2} \left|(f+g)_{xxx}\right|^{p+1}+\frac{s^+}{3\mu_0^+}g^3|(f+g)_{xxx}|^2\,dx\\
&+\int_\Omega f\left|\frac{1}{2\sqrt{ms^+}}f\left(s^- f_{xxx}+ms^+(f+g)_{xxx}\right)+\sqrt{ms^+}g(f+g)_{xxx}\right|^2\,dx\\
&+\int_\Omega\frac{1}{12 ms^+}f^3\left|s^- f_{xxx}  + ms^+ (f+g)_{xxx}\right|^2\,dx,
\end{split}
\end{align}
where we use the definition $\Phi(d)=|d|^{p-1}d$.
In particular, we find that $D(f,g)\geq 0$ and that system \eqref{eq:system_new} is dissipative in the sense that \sout{with}
\begin{align*}
	\frac{d}{dt} E(f,g)(t)
	&=-D(f,g)(t), \quad t \in (0,T).
\end{align*}
Integration with respect to time  yields the desired energy identity \eqref{eq:energy_id}.
% \begin{equation}\label{eq:E}
% E(f,g)(t)+\int_0^tD(f,g)(\tau)\,d\tau = E(f_0,g_0).
% \end{equation}

\end{proof}

In order to prove uniqueness of strong solutions to \eqref{eq:system_new} which  emanate  from the same initial datum let us introduce the relative energy of  two  solutions $(f,g)$ and $(F,G)$ of \eqref{eq:system_new}  by
\begin{equation}\label{eq:rel_energy}
E_{\text{rel}}\bigl((f,g),(F,G)\bigr):=\frac{1}{2}\int_{\Omega}\left|(F+G)_x - (f+g)_x\right|^2+\frac{s^-}{ms^+}\left|(F-f)_x\right|^2\,dx.
\end{equation}

\begin{theorem}[Uniqueness]\label{thm:uniqueness}
Let $p > 1$. Moreover, let $(f,g)$ and $(F,G)$ be two positive solutions of \eqref{eq:system_new} on $[0,T]$, as in Theorem \ref{thm:existence}, emanating from the same initial value $(f_0,g_0)$ with $f_0,g_0 > 0$ on $\bar{\Omega}$. Then $(f,g)=(F,G)$ on $[0,T]$.
\end{theorem} 

\begin{proof} Let $(f,g)$ and $(F,G)$ be two solutions of \eqref{eq:system_new} in the sense of 
Theorem \ref{thm:existence} on some  time interval $[0,T]$, with 
\begin{equation}\label{eq:initial}
    (f_0,g_0)=(F_0,G_0) >0, \quad x \in \bar{\Omega}.
\end{equation}
In view of a continuation argument, it is sufficient to show that there exists a time $T_*\leq T$ such that $(f,g)=(F,G)$ on $[0,T_*]$. 
%Due to Lemma \ref{lem:reg}, the relative energy defined in \eqref{eq:rel_energy} is differentiable with respect to time. 
Following the approach in \cite{LM:2020}, which was applied in the context of a single non-Newtonian thin film, we  show that
\begin{equation}\label{eq:aim}
    \sup_{s\in (0,t)}E_{\text{rel}}\bigl((f,g),(F,G)\bigr)(s)\leq C(t)\sup_{s\in (0,t)}E_{\text{rel}}\bigl((f,g),(F,G)\bigr)(s),
\end{equation}
for some positive constant $C(t)$ depending on $t$ with $C(t)<1$ for $t$ small enough. Note that then 
\[\sup_{s\in (0,t)}E_{\text{rel}}\bigl((f,g),(F,G)\bigr)(s)=0\]
 for $t$ small enough, which in turn implies that $F-f$ and $G-g$ are constant on $(0,t)$. Since $(f,g)$ and $(F,G)$ are continuous with respect to time and have  the same initial datum, we conclude that $(f,g)=(F,G)$ on $(0,t)$, which proves the statement. 
 
In the remainder of the proof we show that \eqref{eq:aim} holds true.
%  The remainder of the proof is devoted to that \eqref{eq:aim} holds. 
 We write
 \begin{align*}
 E_{\text{rel}}&\bigl((f,g),(F,G)\bigr)(t)=\frac{1}{2}\int_{\Omega}\left|(F+G)_x - (f+g)_x\right|^2+\frac{s^-}{ms^+}\left|(F-f)_x\right|^2\, dx.\\
 &=E(F,G)(t)+E(f,g)(t)- \int_\Omega (F+G)_{x}(f+g)_x+\frac{s^-}{ms^+}F_x f_x\, dx \\
 &=-\int_0^t D(F,G)(\tau)\,d\tau -\int_0^t D(f,g)(\tau)\,d\tau-\int_0^t \frac{d}{dt}\int_\Omega (F+G)_{x}(f+g)_x+\frac{s^-}{ms^+}F_xf_x\,dx\,d\tau.
 \end{align*}
 where we use \eqref{eq:energy_id} and \eqref{eq:initial}. Observe that
 \begin{align*}
%  \frac{d}{dt}&\int_\Omega (F+G)_{x}(f+g)_x+\frac{s^-}{ms^+}F_xf_x\,dx=-\int_\Omega (F+G)_{t}(f+g)_{xx}+(F+G)_{xx}(f+g)_{t}\,dx \\
%  &\hspace{6cm}+\frac{s^-}{ms^+}\int_{\Omega}F_{t}f_{xx}+F_{xx}f_{x}\,dx\\
-\,\frac{d}{dt}&\int_\Omega (F+G)_{x}(f+g)_x+\frac{s^-}{ms^+}F_xf_x\, dx \\
 &=\int_\Omega (F+G)_{t}(f+g)_{xx}+(F+G)_{xx}(f+g)_{t}\,dx
 +\frac{s^-}{ms^+}\int_{\Omega}F_{t}f_{xx}+F_{xx}f_{t}\,dx\\
 &=\int_{\Omega}\bigl( p_{21}(F,G)F_{xxx} + p_{22}(F,G)(F+G)_{xxx} + C_p |G|^{p +2}\Phi((F+G)_{xxx}) \bigr)(f+g)_{xxx}\,dx\\
 &+\int_{\Omega} \bigl( p_{21}(f,g)f_{xxx} + p_{22}(f,g)(f+g)_{xxx}  + C_p |g|^{p +2}\Phi((f+g)_{xxx})\bigr)(F+G)_{xxx}\,dx\\
 &+\frac{s^-}{ms^+}\int_{\Omega}\bigl(  p_{11}(F,G)F_{xxx} + p_{12}(F,G)(F+G)_{xxx}\bigr)f_{xxx}\,dx\\
 &+\frac{s^-}{ms^+}\int_{\Omega}\bigl(  p_{11}(f,g)f_{xxx} + p_{12}(f,g)(f+g)_{xxx}\bigr)F_{xxx}\,dx,
 \end{align*}
 where we use that $(f,g)$ and $(F,G)$ solve \eqref{eq:system_new} and satisfy the Neumann-type boundary conditions incorporated in $H^4_B(\Omega)\times H^4_B(\Omega)$.  In the following  we use the notation $\int_{\Omega_t}:=\int_0^t \int_{\Omega}$. Recalling the definition of the dissipative term $D$ in Lemma \ref{lem:Energy}, we find that
 \begin{align}\label{eq:Erel_estimate}
 \begin{split}
  E_{\text{rel}}&\bigl((f,g),(F,G)\bigr)(t)\\
   &=-\int_{\Omega_t}\left(p_{21}(F,G)F_{xxx}-p_{21}(f,g)f_{xxx}\right) \left((F+G)_{xxx}-(f+g)_{xxx}\right)\,d(x,\tau)\\
     &-\int_{\Omega_t}\left(p_{22}(F,G)(F+G)_{xxx}-p_{22}(f,g)(f+g)_{xxx}\right) \left((F+G)_{xxx}-(f+g)_{xxx}\right)\,d(x,\tau)\\      
   &-\frac{s^-}{ms^+}\int_{\Omega_t}\left(p_{11}(F,G)F_{xxx}-p_{11}(f,g)f_{xxx}\right) \left(F_{xxx}-f_{xxx}\right)\,d(x,\tau)\\
   &-\frac{s^-}{ms^+}\int_{\Omega_t}\left(p_{12}(F,G)(F+G)_{xxx}-p_{12}(f,g)(f+g)_{xxx}\right) \left(F_{xxx}-f_{xxx}\right)\,d(x,\tau)\\
   &-C_p\int_{\Omega_t}\left(|G|^{p +2}\Phi((F+G)_{xxx})-|g|^{p +2}\Phi((f+g)_{xxx})\right) \left((F+G)_{xxx}-(f+g)_{xxx}\right)\,d(x,\tau).
   \end{split}
 \end{align}
In order to lighten the notation we set $A:=F_{xxx}-f_{xxx}$ and $B:= G_{xxx}-g_{xxx}$. 
We rewrite each of the five integrals in \eqref{eq:Erel_estimate}  in the same spirit as  the first integral:
\begin{align*}
\int_{\Omega_t}&\left(p_{21}(F,G)F_{xxx}-p_{21}(f,g)f_{xxx}\right) \left(A+B\right)\,d(x,\tau)\\
&=\int_{\Omega_t}p_{21}(F,G)A \left(A+B \right)\,d(x,\tau)+\int_{\Omega_t}\left(p_{21}(F,G)-p_{21}(f,g)\right)f_{xxx} \left(A+B\right)\,d(x,\tau).
\end{align*}
Special attention should be paid to the last integral, which  becomes
\begin{align*}
&\hspace{-1cm}\int_{\Omega_t}\left(|G|^{p +2}\Phi((F+G)_{xxx})-|g|^{p +2}\Phi((f+g)_{xxx})\right) \left((F+G)_{xxx}-(f+g)_{xxx}\right)\,d(x,\tau)\\
&\hspace{-1cm}=\int_{\Omega_t}|G|^{p +2}\left(|(F+G)_{xxx}|^{p-1}(F+G)_{xxx}-|(f+g)_{xxx}|^{p-1}(f+g)_{xxx}\right)\left((F+G)_{xxx}-(f+g)_{xxx}\right)\,d(x,\tau)\\
&\hspace{-1cm}+\int_{\Omega_t}\left(|G|^{p +2}-|g|^{p +2}\right)|(f+g)_{xxx}|^{p-1}(f+g)_{xxx}\left((F+G)_{xxx}-(f+g)_{xxx}\right)\,d(x,\tau).
\end{align*}
Due to the  inequality (cf. \cite[Lemma 1.4.4]{DiB})
\begin{align*}
    \left(\left|(F+G)_{xxx}\right|^{p-1}(F+G)_{xxx}\right. & \left. -\left|(f+g)_{xxx}\right|^{p-1}(f+g)_{xxx}\right)\left((F+G)_{xxx}-(f+g)_{xxx}\right)\\
    &\quad \geq c_p \left|(F+G)_{xxx}-(f+g)_{xxx}\right|^{p+1},
\end{align*}
for some $c_p>0$,
and since $G$ is bounded away from zero, we may estimate
\begin{align*}
    \int_{\Omega_t}&\left(|G|^{p +2}\Phi\bigl((F+G)_{xxx}\bigr)-|g|^{p +2}\Phi\bigl((f+g)_{xxx}\bigr)\right) \left((F+G)_{xxx}-(f+g)_{xxx}\right)\,d(x,\tau)\\
    &\geq c_p\int_{\Omega_t}|G|^{p +2}\left|(F+G)_{xxx}-(f+g)_{xxx}\right|^{p+1}\,d(x,\tau)\\
    &+\int_{\Omega_t}\left(|G|^{p +2}-|g|^{p +2}\right)|(f+g)_{xxx}|^{p-1}(f+g)_{xxx}\left((F+G)_{xxx}-(f+g)_{xxx}\right)\,d(x,\tau).
\end{align*}
Rearranging the terms in \eqref{eq:Erel_estimate}, we obtain
\begin{align*}
    \hspace{-1cm}E_{\text{rel}}&\bigl((f,g),(F,G)\bigr)(t)\\
    &\hspace{-1cm}\leq-\int_{\Omega_t}p_{22}(F,G)(A+B)^2+ \left(p_{21}(F,G)+\frac{s^-}{ms^+}p_{12}(F,G)\right)A(A+B)+\frac{s^-}{ms^+}p_{11}(F,G)A^2\,d(x,\tau)\\
    &\hspace{-1cm}-\int_{\Omega_t}\left(p_{21}(F,G)-p_{21}(f,g)\right)f_{xxx} (A+B)\,d(x,\tau)-\int_{\Omega_t}\left(p_{22}(F,G)-p_{22}(f,g)\right)(f+g)_{xxx}(A+B)\,d(x,\tau)\\
   &\hspace{-1cm}-\frac{s^-}{ms^+}\int_{\Omega_t} \left(p_{11}(F,G)-p_{11}(f,g)\right)f_{xxx}A\,d(x,\tau)-\frac{s^-}{ms^+}\int_{\Omega_t}\left(p_{12}(F,G)-p_{12}(f,g)\right)(f+g)_{xxx}A\,d(x,\tau)\\
   &\hspace{-1cm}-C_p c_p\int_{\Omega_t}|G|^{p +2} (A+B)^{p+1}\,d(x,\tau)-C_p\int_{\Omega_t}\left(|G|^{p +2}-|g|^{p +2} \right)|(f+g)_{xxx}|^{p-1}(f+g)_{xxx}(A+B)\,d(x,\tau).
\end{align*}
Note that the integrand of the first integral  is positive and given by the sum of squares as in Lemma \ref{lem:A}. Thus,
\begin{align*}
    & \hspace{-1cm}E_{\text{rel}}\bigl((f,g),(F,G)\bigr)(t)+ \int_{\Omega_t}f\left(\frac{1}{2\sqrt{ms^+}}f\left(s^- A+ms^+(A+B)\right)+\sqrt{ms^+}g(A+B)\right)^2\,d(x,\tau)\\
    &\hspace{-1cm}+\int_{\Omega_t}\frac{1}{12 ms^+}f^3\left(s^- A  + ms^+ (A+B)\right)^2+\frac{s^+}{3\mu_0^+}g^3(A+B)^2+C_p |G|^{p +2} (A+B)^{p+1}\,d(x,\tau)\\
    &\hspace{-1cm}\leq-\int_{\Omega_t}\left(p_{21}(F,G)-p_{21}(f,g)\right)f_{xxx} (A+B)\,d(x,\tau)-\int_{\Omega_t}\left(p_{22}(F,G)-p_{22}(f,g)\right)(f+g)_{xxx}(A+B)\,d(x,\tau)\\
    &\hspace{-1cm}-\frac{s^-}{ms^+}\int_{\Omega_t} \left(p_{11}(F,G)-p_{11}(f,g)\right)f_{xxx}A\,d(x,\tau)-\frac{s^-}{ms^+}\int_{\Omega_t}\left(p_{12}(F,G)-p_{12}(f,g)\right)(f+g)_{xxx}A\,d(x,\tau)\\
    &\hspace{-1cm}-C_p\int_{\Omega_t}\left(|G|^{p +2}-|g|^{p +2} \right)|(f+g)_{xxx}|^{p-1}(f+g)_{xxx}(A+B)\,d(x,\tau).
 \end{align*}
% \begin{align*}
% p_{22}&(f,g)(A+B)^2+ \left(p_{21}(f,g)+\frac{s^-}{ms^+}p_{12}(f,g)\right)A(A+B)+\frac{s^-}{ms^+}p_{11}(f,g)A^2\\
% &=\frac{\sigma^+}{3\mu^+m^+}g^3(A+B)^2 +\frac{1}{12 m^-\mu^- \sigma^+}f^3\left(\sigma^- A+\sigma^+(A+B)\right)^2\\
% &+\frac{1}{m^-\mu^-}f\left(\frac{1}{2\sqrt{\sigma^+}}f\left(\sigma^- A+\sigma^+(A+B)\right)+\sqrt{\sigma^+}g(A+B)\right)^2.
% \end{align*}
 By elementary estimates, see Lemma \ref{lem:B}, there exists a constant $C>0$ such that
\begin{equation}\label{eq:aux_1}
    \left|p_{ij}(F,G)-p_{ij}(f,g)\right|^2\leq C \bigl((F+G)-(f+g)\bigr)^2 %\quad \mbox{for}\qquad 1\leq i,j\leq 2.
\end{equation}
for $1\leq i,j\leq 2$, and
\begin{equation} \label{eq:aux_2}
    \left||G|^{p +2}-|g|^{p +2} \right|^2\leq C \left|G-g\right|^2.
\end{equation}
In view of the inequality $|G-g|\leq|(F+G)-(f+g)|+|F-f|$, we may deduce from \eqref{eq:aux_1} and \eqref{eq:aux_2}  the estimate
 \begin{equation}\label{eq:estimate_max}
 \max \left(\left|p_{ij}(F,G)-p_{ij}(f,g)\right|^2,\left||G|^{p +2}-|g|^{p +2} \right|^2\right)\leq c \left(\bigl((F+G)-(f+g)\bigr)^2+\frac{s^-}{ms^+}(F-f)^2\right),
 \end{equation}
where  $c>0$ is a generic constant, depending on the supremum  norm of $F,G,f,g$. 
Applying Young's inequality and \eqref{eq:estimate_max} to all integrals on the right-hand side of the estimate for the relative energy $E_{\text{rel}}$,
% to all of the above integrals on the right-hand side, 
we deduce that 
\begin{align*}
    & E_{\text{rel}}\bigl((f,g),(F,G)\bigr)(t)+\int_{\Omega_t}f\left(\frac{1}{2\sqrt{ms^+}}f\left(s^- A+ms^+(A+B)\right)+\sqrt{ms^+}g(A+B)\right)^2\,d(x,\tau)\\
    &+\int_{\Omega_t}\frac{1}{12 ms^+}f^3\left(s^- A  + ms^+ (A+B)\right)^2+\frac{s^+}{3\mu_0^+}g^3(A+B)^2+C_p |G|^{p +2} (A+B)^{p+1}\,d(x,\tau)\\
    &\leq\eps \int_{\Omega_t} A^2+(A+B)^2\,d(x,\tau)\\
    &+ \frac{c}{4\eps}\int_{\Omega_t} \left(\bigl((F+G)-(f+g)\bigr)^2+\frac{s^-}{ms^+}(F-f)^2\right)\left(|(f+g)_{xxx}|^2+|f_{xxx}|^2+|(f+g)_{xxx}|^{2p} \right)\,d(x,\tau) 
\end{align*}
for any $\eps>0$.
As shown in Lemma \ref{lem:C}, there exists $\eps>0$ depending on the supremum  norm of $f$ and $g$ such that 
\begin{align*}
    \eps\left( A^2+(A+B)^2\right)&\leq \frac{1}{12 ms^+}f^3\left(s^- A  + ms^+ (A+B)\right)^2 + \frac{s^+}{3\mu_0^+}g^3(A+B)^2.
\end{align*}
Fixing such an $\eps > 0$, we find that
\begin{align*}
    &\sup_{s\in (0,t)} E_{\text{rel}}\bigl((f,g),(F,G)\bigr)(t)\\
    &\leq 
    \frac{c}{16\eps} \left\|\bigl((F+G)-(f+g)\bigr)^2+\frac{s^-}{ms^+}(F-f)^2 \right\|_{L_\infty(\Omega_t)} \int_{\Omega_t}\left|(f+g)_{xxx}\right|^2+\left|f_{xxx}\right|^2+\left|(f+g)_{xxx}\right|^{2p}\,d(x,\tau).
\end{align*}
In virtue of the Sobolev embedding $H^1(\Omega)\subset L^\infty(\Omega)$ we eventually arrive at
\begin{align*}
    \sup_{s\in (0,t)}E_{\text{rel}}\bigl((F,G),(f,g)\bigr)(s)
    &\leq 
    C(t) \sup_{s\in (0,t)} E_{\text{rel}}\bigl((F,G),(f,g)\bigr)(s),
\end{align*}
where 
\[
    C(t)=\frac{c}{16\eps}\int_{\Omega_t}\left|(f+g)_{xxx}\right|^2+\left|f_{xxx}\right|^2+\left|(f+g)_{xxx}\right|^{2p}\,d(x,\tau)
\]
is a constant depending continuously on $t$ and satisfying $\lim_{t\to 0}C(t)=0$, since $(f,g)$ is a strong solution.
Choosing $t$ small enough, we conclude that $(F_x,G_x)=(f_x,g_x)$ which yields the assertion.
\end{proof}

%=============================================================================
%=============================================================================
%=============================================================================
\bigskip

\section{Long-time behaviour}\label{sec:long-time}

In the previous part of the paper we have proved that, given general positive initial film heights, there exists a unique strong solution to the two-phase thin-film problem \eqref{eq:system} on a small time interval $(0,T)$. In this section we are concerned with the behaviour of solutions for long times. More precisely, we first characterise the maximal time $\bar{T}$ of existence of solutions by showing that, when solutions cease to exists, they either blow up in some intermediate norm or one of the film heights drops down to zero. Second, we prove that the set of steady-state solutions of \eqref{eq:system} consists of positive constants only, i.e. only flat films are steady states of \eqref{eq:system}. Finally, for flow behaviour exponents $p \geq 2$ we show that these flat films are stable in the sense that solutions which are initially close to a positive constant, converge to this constant as time tends to infinity. The restriction to flow behaviour exponents $p \geq 2$ is again due to the previously mentioned lack of regularity. More precisely, the study of stable and center manifolds is usually based on a contraction argument for a related fixed-point problem. However, possible remedies are beyond the scope of this paper.
%=============================================================================
\bigskip

\subsection{Maximal time of existence}

In this part we prove that solutions with finite lifetime only do either blow up in some $H^{4\gamma}_B(\Omega)$-norm or they develop a film rupture at some point $x \in \bar{\Omega}$. Note that we do not rule out the possibility that both scenarios happen simultaneously. To be precise, let $\theta > \rho > 7/8$.
For any initial value $u_0 = (f_0,g_0) \in H^{4\theta}_B(\Omega)\times H^{4\theta}_B(\Omega)$ with $f_0(x),g_0(x) > 0$ for all $x \in \bar{\Omega}$, we first define the maximal time $\bar{T}$ of existence of solutions to \eqref{eq:system} by
\begin{equation*}
	\bar{T} = \sup\left\{T > 0;\ \text{there exists a solution $u=(f,g)$ of \eqref{eq:system} in the sense of Theorem \ref{thm:existence}}\right\}.
\end{equation*}
The uniqueness result, Theorem \ref{thm:uniqueness}, implies that there exists a solution 
\begin{equation*}
	u=(f,g) \in C\bigl([0,\bar{T});H^{4\theta}_B(\Omega)\times H^{4\theta}_B(\Omega)\bigr) \cap C\bigl((0,\bar{T});H^4_B(\Omega)\times H^4_B(\Omega)\bigr) \cap C^1\bigl((0,\bar{T});L_2(\Omega)\times L_2(\Omega)\bigr).
\end{equation*}
We are now prepared to prove the following result.

\begin{theorem}
Let $u_0 = (f_0,g_0) \in H^{4\theta}_B(\Omega)\times H^{4\theta}_B(\Omega)$ with $f_0(x),g_0(x) > 0$ for all $x \in \bar{\Omega}$. Suppose that $\bar{T} < \infty$. Then
\begin{equation}\label{eq:blow-up}
	\liminf_{\bar{T} \nearrow \infty} \frac{1}{\min_{\bar{\Omega}} \{f(t), g(t)\}} + \left\|(f(t),g(t))\right\|_{H^{4\gamma}_B(\Omega)} = \infty
\end{equation}
for all $\gamma \in (\rho,1]$.
\end{theorem}

The proof is based on the standard continuation argument. We follow the lines of the proof of \cite[Thm. 7.1]{LM:2020}.

\begin{proof}
 We fix $\gamma \in (\rho,\theta]$ and assume by contradiction that \eqref{eq:blow-up} is false. Consequently, there exist a sequence $(t_k)_{k\in \N} \subset \R_+$ with $\lim_{k\to \infty} t_k = \bar{T}$ and constants $m, M = M(\gamma) > 0$ such that
\begin{equation*}
	\min_{x \in \bar{\Omega}} \{f(t_k),g(t_k)\} > m \quad \text{and} \quad \left\|(f(t_k),g(t_k))\right\|_{H^{4\gamma}_B(\Omega)} < M, \quad k \in \N.
\end{equation*}
Since $\gamma \leq \theta$ and $\nu = \frac{\gamma-\rho}{2} < \theta-\rho$, we may invoke
Theorem \ref{thm:existence} to deduce that there exists a positive time $T=T(m,M) > 0$ such that, for all $k\in\N$ and corresponding initial values $U_k(0) =(F_k(0),G_k(0))=(f(t_k),g(t_k))$, there exists a solution $U_k = (F_k,G_k)$, satisfying
\begin{equation*}
	U_k \in C\bigl([0,T];H^{4\gamma}_B(\Omega)\times H^{4\gamma}_B(\Omega)\bigr) \cap C^\frac{\gamma-\rho}{2}\bigl([0,T];H^{4\rho}_B(\Omega)\times H^{4\rho}_B(\Omega)\bigr) \cap C\bigl((0,T];H^4_B(\Omega)\times H^4_B(\Omega)\bigr). 
\end{equation*}
This implies that
\begin{equation*}
	\begin{cases}
		\Acal\circ U_k \in C^\alpha\bigl([0,T];\Hcal(H^4_B(\Omega)\times H^4_B(\Omega);L_2(\Omega)\times L_2(\Omega))\bigr) & \\ 
		\Fcal\circ U_k \in C^\alpha\bigl([0,T];L_2(\Omega)\times L_2(\Omega)\bigr), &
	\end{cases}
\end{equation*}
where $\alpha = (p-1)\frac{\gamma-\rho}{2}$. Moreover, we have that $U_k(0)=(f(t_k),g(t_k)) \in H^4_B(\Omega)$. Due to this property,  linear parabolic regularity yields 
\begin{equation*}
	U_k \in C\bigl([0,T];H^4_B(\Omega)\times H^4_B(\Omega)\bigr) \cap C^1\bigl([0,T];L_2(\Omega)\times L_2(\Omega)\bigr).
\end{equation*}
This allows us to extend the solution at time $t_k$ as follows. For $t_k \geq \bar{T} - \frac{T}{2}$ we define the extension $\tilde{u} = (\tilde{f},\tilde{g})$ by
\begin{equation*}
	\tilde{u} = \begin{cases}
	u(t) = (f(t),g(t)), & t \in [0,t_k) \\
	U_k(t - t_k), & t \in [t_k,t_k+T].
	\end{cases}
\end{equation*}
Then, the extended function
\begin{equation} \label{eq:reg_extension}
	\tilde{u} \in C\bigl([0,t_k+T];H^{4\theta}_B(\Omega)\times H^{4\theta}_B(\Omega)\bigr) \cap C\bigl((0,t_k+T];H^4_B(\Omega)\times H^4_B(\Omega)\bigr) \cap C^1\bigl((0,t_k+T];L_2(\Omega)\times L_2(\Omega)\bigr)
\end{equation}
is a solution to the problem
\begin{equation} \label{eq:ext_prob}
	\begin{cases}
	\tilde{u}_t + \Acal(\tilde{u}) \tilde{u} = \Fcal(\tilde{u}), & t > 0, \\
	\tilde{u}(0) = u(t_n) &
	\end{cases}
\end{equation}
on $(0,t_k+T)$. Indeed, clearly, $\tilde{u}$ is continuous in the sense that $\tilde{u} \in C\bigl([0,t_k+T];H^{4\theta}_B(\Omega)\times H^{4\theta}_B(\Omega)\bigr) \cap C\bigl((0,t_k+T];H^4_B(\Omega)\times H^4_B(\Omega)\bigr)$ and the differential equation
\begin{equation*}
	\tilde{u}_t + \Acal(\tilde{u}) \tilde{u} = \Fcal(\tilde{u})
\end{equation*}
is satisfied on both intervals $(0,t_k)$ and $(t_k,t_k+T)$. Due to the regularity $\tilde{u} \in C\bigl((0,t_k+T];H^4_B(\Omega)\times H^4_B(\Omega)\bigr)$, we have in addition that the limits 
\begin{equation*}
	\begin{split}
		0 &= \lim_{t \searrow t_k} \partial_t^+ U_k(t-t_k) + \Acal\bigl(U_k(t-t_k)\bigr) U_k(t-t_k) - \Fcal\bigl(U_k(t-t_k)\bigr)
		\\
		&= \partial_t^+ U_k(0) + \Acal\bigl(U_k(0)\bigr) U_k(0) - \Fcal\bigl(U_k(0)\bigr) \\
		&=
		\partial_t^+ \tilde{u}(t_k) + \Acal\bigl(\tilde{u}(t_k)\bigr) \tilde{u}(t_k) - \Fcal\bigl(\tilde{u}(t_k)\bigr) \quad \text{in } L_2(\Omega)\times L_2(\Omega)
	\end{split}
\end{equation*}
for $t > t_k$, and
\begin{equation*}
	\begin{split}
	0 &= \lim_{t \nearrow t_k} \partial_t^- u(t) + \Acal(u(t)) u(t) - \Fcal(u(t)) \\
	&=
	\partial_t^- \tilde{u}(t_k) + \Acal\bigl(\tilde{u}(t_k)\bigr) \tilde{u}(t_k) - \Fcal\bigl(\tilde{u}(t_k)\bigr) \quad \text{in } L_2(\Omega)\times L_2(\Omega)
	\end{split}
\end{equation*}
for $t < t_k$ are well-defined. Therefore, $\tilde{u}_t$ can be uniquely extended at time $t_k$ to a function $C\bigl([0,t_k+T];L_2(\Omega)\times L_2(\Omega)\bigr)$. This proves that $\tilde{u}$ is a solution of the extended problem \eqref{eq:ext_prob} on $(0,t_k+T)$ in the sense of \eqref{eq:reg_extension}. Since $t_k+T \geq \bar{T}+\frac{T}{2}$, this is a contradiction to the maximality of $\bar{T}$. 
\end{proof}

%=============================================================================
\bigskip

\subsection{Characterisation of steady states}

This section is concerned with stability properties of steady state solutions to \eqref{eq:system}. First, we show that flat films of positive height are the only possible steady states. Moreover, using linearised stability, we can prove that solutions, which are initially close to a flat film, have an infinite lifetime and converge to the flat film as time tends to infinity.

%-----------------------

\begin{lemma}[Conservation of mass]\label{lem:CoM}
Let $(f,g)$ be a solution to \eqref{eq:system} $[0,T]$ corresponding to the initial datum $(f_0,g_0)$, as found in Theorem \ref{thm:existence}. Then the solution conserves its mass in the sense that
\[
\int_{\Omega}f(t,x)\,dx = \int_{\Omega}f_0(x)\,dx\qquad \mbox{and}\qquad \int_{\Omega}g(t,x)\,dx = \int_{\Omega}g_0(x)\,dx
\]
for all $t\in [0,T]$.
\end{lemma}

%-----------------------

\begin{proof}
The statement follows immediately from the structure of the equations in \eqref{eq:system} and the Neumann-type boundary conditions \eqref{eq:boundary}.
\end{proof}

%-----------------------

\begin{corollary} \label{cor:steady_states}
A solution $(f,g)$ of \eqref{eq:system} as found in Theorem \ref{thm:existence} is a steady state if and only if $f$ and $g$ are constant. That is, the positive flat two-phase thin film is the only equilibrium solution of \eqref{eq:system} with boundary conditions \eqref{eq:boundary}. 
\end{corollary}

%-----------------------

\begin{proof}
It is clear from the structure of the equation \eqref{eq:system} and the boundary conditions \eqref{eq:boundary} that any pair of positive constants $(f,g)$ is an equilibrium solution. We now show that any steady-state solution of \eqref{eq:system} with boundary conditions \eqref{eq:boundary} 
%is determined to be constant 
must automatically be a positive constant. Recall from Lemma \ref{lem:Energy} that the energy functional satisfies
\begin{equation*}
\frac{d}{dt}E(f,g)(t)=-D(f,g)(t),\quad t \in (0,T),
\end{equation*}
where $D(f,g)$ can be written as 
\begin{equation} \label{eq:diss_steady_states}
	\begin{split}
	D(f,g)&=\int_\Omega C_p g^{p+2} \left|(f+g)_{xxx}\right|^{p+1}+\frac{s^+}{3\mu_0^+}g^3|(f+g)_{xxx}|^2\,dx\\
	&\quad+\int_\Omega f\left|\frac{1}{2\sqrt{ms^+}}f\left(s^- f_{xxx}+ms^+(f+g)_{xxx}\right)+\sqrt{ms^+}g(f+g)_{xxx}\right|^2\,dx\\
	&\quad+\int_\Omega\frac{1}{12 ms^+}f^3\left|s^- f_{xxx}  + ms^+ (f+g)_{xxx}\right|^2\,dx.
	\end{split}
\end{equation}
Note that all three integrands in \eqref{eq:diss_steady_states} are nonnegative. If $(f,g)$ is an equilibrium solution, then $\frac{d}{dt}E(f,g)(t)=0$ for all times and hence $D(f,g)(t)=0$ for all times. But this is only possible if all three integrands 
%in $D(f,g)$ 
on the right-hand side of \eqref{eq:diss_steady_states}
vanish. Since $f,g>0$, we infer in particular that
\[
(f+g)_{xxx}=0\quad \mbox{and}\quad (s^-f_{xxx}+ms^+(f+g)_{xxx})=0.
\]
Together this implies that $f_{xxx}=g_{xxx}=0$ and thus, $f_{xx}, g_{xx}$ are constants. This in turn implies that $f_x$ and $g_x$ are linear. Due to the Neumann boundary conditions, we deduce that $f_x=g_x=0$ and finally, $(f,g)$ is a pair of positive constants.
%In view of Lemma \ref{lem:reg}, the derivative $\frac{d}{dt} E(f,g)$ is well-defined for the solution $(f,g)$ of \eqref{eq:system} and we have
%\begin{equation*}
%	\frac{d}{dt} E(f,g) = - \int_{\Omega} \textcolor{magenta}{diss. term}\, dx = 0.
%\end{equation*}
%Since $u_\ast=(f_\ast,g_\ast) > 0$, we deduce that $(u_\ast)_{xxx} = 0$.
% This implies that the second derivative $(u_\ast)_{xx}$ is constant and consequently, the first derivative $(u_\ast)_{x}$ is linear. Matching with the Neumann boundary condition yields $(u_\ast)_{x}\equiv 0$ and finally, $u_\ast\equiv 0$.
\end{proof}

%-----------------------

\subsection{Behaviour close to stationary solutions} 
In this section we apply the version \cite[Theorem 2.20]{HI} of the center manifold theorem established in \cite{Milke}
for quasilinear problems in Hilbert spaces in order to prove that flat films of positive height are exponentially stable in the case $p\ge 2$. To this end, we eliminate the central spectrum of the linear part of the differential operator by exploiting the conservation of mass. Note that we assume $p\ge 2$ to ensure sufficient regularity of the nonlinear part. We prove the following result.

\begin{theorem} \label{Thm:Asympt_stab}
Let  $p\ge 2$ and let $u_{\ast} \in \R_{+}^2$. Then there exist positive constants $\eps, M, \kappa > 0$ such that for all initial values $u_0 \in H^4_B(\Omega)\times H^4_B(\Omega)$, satisfying 
\begin{equation*}
	u_0(x) > 0 \quad \text{for} \quad x \in \bar{\Omega},
	\quad
	\frac{1}{|\Omega|} \int_{\Omega} u_0(x)\, dx = u_{\ast}
	\quad
	\text{and}
	\quad
	\left\|u_0 - u_{\ast}\right\|_{H^4_B(\Omega)\times H^4_B(\Omega)} \leq \eps,
\end{equation*}
the solution $u$ of \eqref{eq:system} exists globally in time and we have the exponential bound
\begin{equation*}
	\left\|u(t) - u_{\ast}\right\|_{H^4_B(\Omega)\times H^4_B(\Omega)} + \left\|u_t(t)\right\|_{L_2(\Omega)\times L_2(\Omega)}
	\leq
	M \exp^{-\kappa t} \left\|u_0 - u_{\ast}\right\|_{H^4_B(\Omega)\times H^4_B(\Omega)}.
\end{equation*}
\end{theorem}
To prove Theorem \ref{Thm:Asympt_stab}, we need to verify that Problem \eqref{eq:system} fits into the abstract setting of \cite[Theorem 2.20]{HI}. For this purpose, we start by introducing some notations and requirements. 

\medskip

\subsection*{Projection of the original problem} 
From Lemma \ref{lem:CoM} we know that each solution $u$ of \eqref{eq:Cauchy_Last} conserves its mass. According to that, we introduce the notation
\begin{equation*}
	u_\ast = \frac{1}{|\Omega|}\int_{\Omega} u(x)\, dx = \frac{1}{|\Omega|}\int_{\Omega} u_0(x)\, dx  > 0
\end{equation*}
for the average film heights $u_\ast=(f_\ast,g_\ast) \in \R^2_+$ of $u$ and incorporate this property into our analysis by defining the projection
%We now incorporate the mass preserving property, cf. Lemma \ref{lem:CoM}, into the asymptotic analysis. That is, we eliminate the non-zero constant functions from the spaces we work in. To this end, we introduce the projection
\begin{equation*}
\begin{cases}
	P \in \Lcal\bigl(L_2(\Omega)\times L_2(\Omega)\bigr) \cap \Lcal\bigl(H^4_B(\Omega)\times H^4_B(\Omega)\bigr) & \\	
	P u := u - u_\ast.
\end{cases}
\end{equation*}
Via $P$ we can decompose $L_2(\Omega)\times L_2(\Omega)$ and $H^4_B(\Omega)\times H^4_B(\Omega)$ into 
direct sums
\begin{equation*}
\begin{cases}
	L_2(\Omega)\times L_2(\Omega) = P \bigl(L_2(\Omega)\times L_2(\Omega)\bigr) \oplus (1-P) \bigl(L_2(\Omega)\times L_2(\Omega)\bigr) & \\
	H^4_B(\Omega)\times H^4_B(\Omega) = P \bigl(H^4_B(\Omega)\times H^4_B(\Omega)\bigr) \oplus (1-P) \bigl(H^4_B(\Omega)\times H^4_B(\Omega)\bigr). &
\end{cases}
\end{equation*}
Observe that 
\begin{equation*}
	\begin{cases}
	L_{2,a}(\Omega)\times L_{2,a}(\Omega) := \left\{v \in L_{2}(\Omega)\times L_{2}(\Omega);\, v_\ast=0\right\} =
	P\bigl(L_2(\Omega)\times L_2(\Omega)\bigr) \\
	H^4_{B,a}(\Omega)\times H^4_{B,a}(\Omega) := \left\{v \in H^4_B(\Omega)\times H^4_B(\Omega);\, v_\ast=0\right\} =  P\bigl(H^4_B(\Omega)\times H^4_B(\Omega)\bigr)
	\end{cases}
\end{equation*}
contain the non-constant functions and the zero function of 
$L_2(\Omega)\times L_2(\Omega)$ and $H^4_B(\Omega)\times H^4_B(\Omega)$, respectively, while $(1-P)\bigl(L_2(\Omega)\times L_2(\Omega)\bigr)$ and $(1-P)\bigl(H^4_B(\Omega)\times H^4_B(\Omega)\bigr)$ contain the constant functions of 
$L_2(\Omega)\times L_2(\Omega)$ and $H^4_B(\Omega)\times H^4_B(\Omega)$, respectively.

Given initial film heights $u_0=(f_0,g_0) \in H^{4}_B(\Omega)\times H^{4}_B(\Omega)$, with $u_0(x) > 0, x \in \bar{\Omega}$, let $u=(f,g)$ be the unique strong solution to the Cauchy problem
\begin{equation}\tag{CP}\label{eq:Cauchy_Last}
	\begin{cases}
		u_t + \Acal(u) u = \Fcal(u), \quad t > 0 & \\
		u(0) = u_0,
	\end{cases}
\end{equation}
as obtained by Theorem \ref{thm:existence} and Theorem \ref{thm:uniqueness}. Recall from the proof of Theorem \ref{thm:existence} that, for $u \in H^{4\rho}_B(\Omega)\times H^{4\rho}_B(\Omega)$ with $u(x) \geq \kappa > 0$ componentwise for all $x \in \bar{\Omega}$, the differential operator $\Acal$ satisfies
\begin{equation*}
	\begin{cases}
	\Acal(u) \in \Hcal\left(H^4_B(\Omega)\times H^4_B(\Omega);L_2(\Omega)\times L_2(\Omega)\right) & \\
	\Acal(u) v = A(u,u_{xxx}) \partial_x^4 v
	\end{cases}
\end{equation*}
and the mapping $u\mapsto \Acal(u)$ is continuously differentiable in view of
\begin{equation*}
	\Acal \in C^{p-1}_{\text{loc}}\left(H^{4\rho}_B(\Omega)\times H^{4\rho}_B(\Omega);\Hcal\left(H^4_B(\Omega)\times H^4_B(\Omega);L_2(\Omega)\times L_2(\Omega)\right)\right)
\end{equation*}
for all $p \ge 2$.
Moreover, $\Fcal$ satisfies
\begin{equation*}
	\Fcal \in C^p_{\text{loc}}\left(H^{4\rho}_B(\Omega)\times H^{4\rho}_B(\Omega);L_2(\Omega)\times L_2(\Omega)\right).
\end{equation*} 
Since $u$ is a solution to \eqref{eq:Cauchy_Last} and  satisfies $(1-P) u_0 = u_\ast$ initially, it follows that $(1-P) u(t) = u_\ast$ as long as it exists.
This allows for all $t \geq 0$ the decomposition
\begin{equation}\label{eq:decomp}
	u(t) = P u(t) + u_\ast =: \varphi(t) + u_\ast 
%	\in P\bigl(H^4_B(\Omega)\times H^4_B(\Omega)\bigr) \oplus (1-P) \bigl(H^4_B(\Omega)\times H^4_B(\Omega)\bigr).
\end{equation}
with
\begin{equation*}
	\varphi(t) \in H^4_{B,a}(\Omega)\times H^4_{B,a}(\Omega) = P\bigl(H^4_B(\Omega)\times H^4_B(\Omega)\bigr) \quad \text{and} \quad u_\ast \in (1-P) \bigl(H^4_B(\Omega)\times H^4_B(\Omega)\bigr).
\end{equation*}
We find that $\varphi$ solves the equation
\begin{equation} \label{eq:diffeq_phi}
	\varphi_t + \Acal_\ast \varphi =	\Fcal_\ast(\varphi)
\end{equation}
with
\begin{equation*}
    \Acal_\ast = \Acal(u_\ast)
    \quad \text{and} \quad
	\Fcal_\ast(\varphi)
	=
	\Fcal(u_\ast +\varphi) - \Fcal(u_\ast) - \left(\Acal(u_\ast + \varphi) - \Acal(u_\ast)\right) \varphi.
\end{equation*}
%{\color{orange} GB: Notice $\Fcal(u_*)=0$.}
Observe that
\begin{equation}\label{eq:projection}
    (1-P)(\Acal_*\phi)=0\qquad \mbox{and} \qquad (1-P)\mathcal{F}_*(\phi)=0\qquad \mbox{for}\quad \phi \in H^4_{B,a}(\Omega)\times H^4_{B,a}(\Omega),
\end{equation}
hence
\begin{equation}
\begin{cases}
    \Acal_* \in \mathcal{L} \bigl(H^4_{B,a}(\Omega)\times H^4_{B,a}(\Omega);L_{2,a}(\Omega)\times L_{2,a}(\Omega)\bigr),&\\
    \Fcal_* \in C^{p-1}_{\text{loc}}\bigl(H^4_{B,a}(\Omega)\times H^4_{B,a}(\Omega);L_{2,a}(\Omega)\times L_{2,a}(\Omega)\bigr),
\end{cases}
\end{equation}
with $\Fcal_\ast(0)=0$ and $D \Fcal_\ast(0) = 0$.
In order to prove the asymptotic stability of $u_*$ for \eqref{eq:Cauchy_Last}, we are left to show that $\phi=0$ is asymptotically stable for \eqref{eq:diffeq_phi}. We start by verifying that the linear autonomous operator $-\Acal_*$ is the infinitesimal generator of an analytic semigroup, which is exponentially stable.

\begin{lemma}\label{lem:spectrum_A_*}
The operator $\Acal_*$ belongs to $\mathcal{H}\bigl(H^4_{B,a}(\Omega)\times H^4_{B,a}(\Omega);L_{2,a}(\Omega)\times L_{2,a}(\Omega)\bigr)$ and there exists $\kappa>0$ such that
\[
\sigma(-\Acal_*)\subset \{\lambda \in \C \mid \mbox{\textit{Re}}\,\lambda <-\kappa\}.
\]
\end{lemma}

\begin{proof}
We already know that $\Acal_*\in \mathcal{H}\left(H^4_B(\Omega)\times H^4_B(\Omega);
L_2(\Omega)\times L_2(\Omega)\right)$. In view of \eqref{eq:projection}, we find that the restriction of $-\Acal_*\big|_{H^4_{B,a}(\Omega)\times H^4_{B,a}(\Omega)}$ is the infinitesimal generator of an analytic semigroup on $L_{2,a}(\Omega)\times L_{2,a}(\Omega)$.
Since $\Acal_* \bar \phi=0$ for any constant function $\bar \phi$ on $\Omega$, it is clear that zero belongs to the spectrum of $-\Acal_*$ as an operator on $H^4_{B}(\Omega)\times H^4_B(\Omega)$. However, the restriction of $-\Acal_*$ onto the space $H^4_{B,a}(\Omega)\times H^4_{B,a}(\Omega)$ of functions in $H^4_{B}(\Omega)\times H^4_{B}(\Omega)$ with zero average eliminates the central spectrum $\sigma_0(-\Acal_*):=\{\lambda \in \C;\, \text{Re}\, \lambda=0\}$.
Let us consider the homogeneous Cauchy problem
\begin{equation}\label{eq:Cauchy_Last_2}
	\begin{cases}
		\phi_t + \Acal_* \phi = 0, \quad t > 0 & \\
		\phi(0) = \phi_0
	\end{cases}
\end{equation}
on $L_{2,a}(\Omega)\times L_{2,a}(\Omega)$. Since $-\Acal_*$ is the infinitesimal generator of an analytic semigroup on $L_{2,a}(\Omega)\times L_{2,a}(\Omega)$, there exists for every $\phi_0 \in L_{2,a}(\Omega)\times L_{2,a}(\Omega)$ a unique solution 
\begin{equation*}
   \phi \in C\bigl((0,\infty);H^4_{B,a}(\Omega)\times H^4_{B,a}(\Omega)\bigr)\cap C^1\bigl((0,\infty);L_{2,a}(\Omega)\times L_{2,a}(\Omega)\bigr).  
\end{equation*}
Recall that
\begin{equation*}
	\Acal_\ast \phi = A_\ast \partial_x^4 \phi, \quad \phi \in H^4_{B,a}(\Omega)\times H^4_{B,a}(\Omega),
\end{equation*}
where
\begin{equation*}
	A_\ast = 
	\begin{pmatrix}
	\frac{ms^+ + s^-}{3} f_\ast^3 + \frac{ms^+}{2}f_\ast^2 g_\ast &
	ms^+ \left(\frac{f_\ast^3}{3}+\frac{f_\ast^2 g_\ast}{2}\right)
	\\
	\frac{ms^+ + s^-}{2} f_\ast^2 g_\ast + ms^+ f_\ast g_\ast^2 + \frac{s^+}{3\mu_0^+} g_\ast^3&  \frac{m s^+}{2} f_\ast^2 g_\ast + m s^+ f_\ast g_\ast^2 + \frac{s^+}{3\mu_0^+} g_\ast^3
	\end{pmatrix},	
\end{equation*}
and that we proved in Section \ref{sec:existence} that $A_\ast$ has two distinct positive eigenvalues $0<\lambda_-<\lambda_+$. Consequently we can diagonalize $A_\ast$, i.e., there exists an invertible matrix $U$ such that
\begin{equation}\label{eq:diagonal}
    A_\ast = U^{-1} \text{diag}(\lambda_-,\lambda_+) U.
\end{equation}
Setting $\|v\|_U:=\|Uv\|_{L_2(\Omega)\times L_2(\Omega)}$ defines an equivalent norm on $L_{2,a}(\Omega)\times L_{2,a}(\Omega)$ and we will prove that there exists $\kappa>0$ such that
\begin{equation}\label{eq:growth}
\|\phi\|_U\leq e^{-\kappa t}\|\phi_0\|_U.
\end{equation}
This strictly negative growth bound on the analytic semigroup generated by $-\Acal_\ast$ implies in turn that
\[
s(-\Acal_*):=\sup \{\mbox{Re} \lambda \mid \lambda \in \sigma(-\Acal_*)\}\leq -\kappa.
\]
Let $\phi$ be a solution to \eqref{eq:Cauchy_Last_2}. In order to establish \eqref{eq:growth}, we define $\psi=U\phi$ and note that \eqref{eq:Cauchy_Last_2} and \eqref{eq:diagonal} together imply that
\begin{equation*}
	\begin{cases}
		\psi_t + \text{diag}(\lambda_-,\lambda_+)\partial_x^4 \psi = 0, \quad t > 0 & \\
		\psi(0) = U\phi_0.
	\end{cases}
\end{equation*}
Testing this equation with $\psi$, integrating by parts twice and using the Neumann-type boundary conditions leads to
\begin{equation*}
\frac{1}{2}\frac{d}{dt}\|\psi\|_{L_2\times L_2}^2 + \left( \text{diag}(\lambda_-,\lambda_+) \psi_{xx},\psi_{xx} \right)_{L_2\times L_2}=0.
\end{equation*}
Using that $\text{diag}(\lambda_-,\lambda_+)$ is a symmetric, positive definite matrix,  and applying the Poincar\'e inequality twice, we find that
\[
\frac{d}{dt} \|\psi\|_{L_2(\Omega)\times L_2(\Omega)}^2 \leq - 2 \kappa \|\psi\|_{L_2(\Omega)\times L_2(\Omega)}^2
\]
for some $\kappa>0$. Inequality \eqref{eq:growth} now follows from $\|\phi\|_U = \|\psi\|_{L_2(\Omega)\times L_2(\Omega)}$.
\end{proof}

%
%------------------
%

\subsection*{Proof of Theorem \ref{Thm:Asympt_stab}}
We verify the hypothesis of \cite[Theorem 2.20]{HI}.
To this end, we identify the Hilbert spaces
\begin{equation*}
Z = H^4_{B,a}(\Omega)\times H^4_{B,a}(\Omega) \longhookrightarrow   
X = L_{2,a}(\Omega)\times L_{2,a}(\Omega),
\end{equation*}
where the embedding is continuous. Moreover, we identify $L=-\Acal_\ast$, $R=\Fcal_\ast$ and consider the problem
\begin{equation*}
\begin{cases}
    \frac{d\phi}{dt} = L\phi + R(\phi), & t > 0\\
    \phi(0) = u_0 - u_\ast
\end{cases}
\end{equation*}
on $L_{2,a}(\Omega)\times L_{2,a}(\Omega)$, where we have $\phi(0)\in H^4_{B,a}(\Omega)\times H^4_{B,a}(\Omega)$. Applying \cite[Theorem 2.20]{HI} requires that $L$ and $R$ satisfy the following properties \cite[Hypotheses 2.1, 2.4 and Equ. (2.9)]{HI}:
\begin{itemize}
	\item[(H2.1.i)] $L \in \Lcal\bigl(H^4_{B,a}(\Omega)\times H^4_{B,a}(\Omega);L_{2,a}(\Omega)\times L_{2,a}(\Omega)\bigr)$;
	\item[(H2.1.ii)] for some $k \geq 2$ there exists a neighbourhood $V \subset H^4_{B,a}(\Omega)\times H^4_{B,a}(\Omega)$ of zero such that $R \in C^k\bigl(V;L_{2,a}(\Omega)\times L_{2,a}(\Omega)\bigr)$ and satisfies $$R(0)=0 \quad \text{and}\quad DR(0)=0;$$
	\item[(H2.4.i)] there exists a positive constant $\kappa > 0$ such that $$\sigma_{-}(L) \subset \left\{\lambda \in \C;\, \text{Re}\, \lambda \leq -\kappa\right\};$$
	\item[(H2.4.ii)] $\sigma_0(L)$ consists of finitely many eigenvalues with finite algebraic multiplicity;
    \item[(H2.9)] there exist constants $\omega_0,C > 0$ such that for all $\omega \in \R$ with $|\omega| \geq \omega_0$, we have that $\ii \omega$ belongs to the resolvent set of $L$ and
    $$\left\|(\ii \omega I - L)^{-1}\right\|_{\Lcal(X)} \leq \frac{C}{|\omega|}.$$
\end{itemize}
Observe that (H2.1.i), (H2.4.i) and (H2.4.ii) follow from Lemma \ref{lem:spectrum_A_*}. Furthermore, (H2.1.ii) is satisfied with the exception that we only have $R\in C^{p-1}\bigl(H^4_{B,a}(\Omega)\times H^4_{B,a}(\Omega); L_{2,a}(\Omega)\times L_{2,a}(\Omega)\bigr)$, $p\ge2$. This is however sufficient for our purposes because it still holds that
\begin{equation*}
    \frac{\sup_{\|u\|_Z} \|R^\eps(u)\|_X}{\eps} + \sup_{\|u\|_Z} \|DR^\eps(u)\|_X \longrightarrow 0 \quad \text{as} \quad \eps \to 0,
\end{equation*}
where $R^\eps(u) = R(u) \chi_\eps\left(\|u\|_Z\right)$ is a regularised version of $R$ and $\chi_\eps$ is a smooth cutoff function that is equal to one on $[0,\eps]$ and equal to zero on $[2\eps,\infty)$. This allows us to perform the contraction argument for the underlying fixed point problem, cf. \cite{Milke}. Finally, (H2.9) follows from the standard resolvent estimate for generators of analytic semigroups and Lemma \ref{lem:spectrum_A_*}.

\appendix

\section{Auxiliary result}

\begin{lemma}\label{lem:A} Let the polynomials $p_{ij}$, $1\leq i,j\leq 2$ be as in \eqref{eq:system_new} and $A,B\in \R$. Then,
\begin{align*}
p_{22}&(f,g)(A+B)^2+ \left(p_{21}(f,g)+\frac{s^-}{ms^+}p_{12}(f,g)\right)A(A+B)+\frac{s^-}{ms^+}p_{11}(f,g)A^2\\
&=f\left(\frac{1}{2\sqrt{ms^+}}f\left(s^- A+ms^+(A+B)\right)+\sqrt{ms^+}g(A+B)\right)^2+\frac{1}{12 ms^+}f^3\left(s^- A  + ms^+ (A+B)\right)^2\\
&+\frac{s^+}{3\mu_0^+}g^3(A+B)^2
\end{align*}
\begin{proof}
Simply inserting the definition of $p_{ij}$, we have that
\begin{align*}
p_{22}&(f,g)(A+B)^2+ \left(p_{21}(f,g)+\frac{s^-}{ms^+}p_{12}(f,g)\right)A(A+B)+\frac{s^-}{ms^+}p_{11}(f,g)A^2\\
&=\left(\frac{ms^+}{3}f^3+ms^+(fg^2+f^2g)+\frac{s^+}{3\mu_0^+}g^3\right)(A+B)^2+\left(\frac{2s^-}{3}f^3+s^-f^2g\right)A(A+B)+\frac{(s^-)^2 }{3ms^+}f^3A^2
\end{align*}
Using the identities
\begin{equation*}
	\frac{\left(s^-\right)^2}{ms^+}f^3 A^2 + 2 s^-f^3 A (A+B) + ms^+ f^3(A+B)^2
	=
	\frac{f^3}{ms^+}\left(s^- A  + ms^+ (A+B)\right)^2
\end{equation*}
and
\begin{equation*}
	s^-f^2g A(A+B) + ms^+f^2g(A+B)^2
	=
	f^2g(A+B)(s^-A + ms^+(A+B))
\end{equation*}
we obtain by completing the square
\begin{align*}
&p_{22}(f,g)(A+B)^2+ \left(p_{21}(f,g)+\frac{s^-}{ms^+}p_{12}(f,g)\right)A(A+B)+\frac{s^-}{ms^+}p_{11}(f,g)A^2\\
&=\frac{1}{3m s^+}f^3\left(s^- A  + m s^+ (A+B)\right)^2 +f^2g(A+B)(s^-A + ms^+(A+B)) + ms^+fg^2(A+B)^2\\
&+ \frac{s^+}{3\mu_0^+}g^3(A+B)^2 \\
&=f\left(\frac{1}{2\sqrt{ms^+}}f\left(s^- A+ms^+(A+B)\right)+\sqrt{ms^+}g(A+B)\right)^2+\frac{1}{12 ms^+}f^3\left(s^- A  + ms^+ (A+B)\right)^2\\
&+\frac{s^+}{3\mu_0^+}g^3(A+B)^2
\end{align*}
\end{proof}

\begin{lemma}\label{lem:B}
Let the polynomials $p_{ij},$, $1\leq i,j\leq 2$ be as in \eqref{eq:system_new}. Then, there exists a constant $c>0$, depending on the absolute value of $X,Y,x,y$ such that
\[
|p_{i,j}(X,Y)-p_{i,j}(x,y)|^2 \leq c ((X+Y)-(x+y))^2,\qquad X,Y,x,y\in \R.
\]
and similarly
\[
\left||G|^{p+2}-|g|^{p+2}\right|^2\leq c(G-g)^2
\]
\end{lemma}

\begin{proof}
Notice that $p_{ij}(x,y)=\sum_{n=0}^3a_{n,ij}x^ny^{3-n}$, $1\leq i,j\leq 2$ are polynomials of order three. Here $a_{n,ij}\in \R$ are the coefficients. Hence,
\begin{align*}
|p_{i,j}(X,Y)-p_{i,j}(x,y)|^2 &\leq  4\sum_{n=0}^3\left|a_{n,ij}X^nY^{3-n}-a_{n,ij}x^ny^{3-n}\right|^2\\
&= 4\sum_{n=0}^3\left|a_{n,ij}X^n\left(Y^{3-n}-y^{3-n}\right)+a_{n,ij}y^{3-n}\left(X^n-x^n\right)\right|^2\\
&\leq c\sum_{n=1}^3\left(\left(Y^{n}-y^{n}\right)^2+\left(X^n-x^n\right)^2\right),
\end{align*}
where $c=c(a_{n,ij},|X|,|y|)>0$ is a constant depending on the coefficients $a_{n,ij}$ and the absolute value of $X,y$. In view of the elementary inequality
\begin{equation}\label{eq:elementary}
a^r-b^r\leq ra^{r-1}(a-b),\qquad \mbox{for} \quad a\geq b\geq 0,\qquad r\geq 1,
\end{equation}
we deduce that there exists a constant $C=C(a_{n,ij},|X|,|Y|,|x|,|y|)>0$ such that
\[
|p_{i,j}(X,Y)-p_{i,j}(x,y)|^2\leq C\left((X-x)^2+(Y-y)^2\right).
\]
Moreover, in view of \eqref{eq:elementary} we have that
\begin{align*}
\left||G|^{p+2}-|g|^{p+2}\right|^2\leq \left|(p+2)\max\left(|G|^{p+1},|g|^{p+1}\right)(G-g)\right|^2,
\end{align*}
which proves the statement.
\end{proof}
\end{lemma}

\begin{lemma}\label{lem:C}
Let $f,g$ be positive, then there exists $\eps>0$ such that
 \begin{align*}
 \eps\left( A^2+(A+B)^2\right)&\leq \frac{1}{12 ms^+}f^3\left(s^- A  + ms^+ (A+B)\right)^2 + \frac{s^+}{3\mu_0^+}g^3(A+B)^2.
 \end{align*}
 \end{lemma}
 
  \begin{proof}
 Notice that
 \begin{align*}
 \left(s^- A+ms^+(A+B)\right)^2&= |s^-|^2A^2 + 2s^-ms^+ A(A+B)+|ms^+|^2(A+B)^2\\
 &\geq |s^-|^2A^2 +|ms^+|^2(A+B)^2-\frac{1}{1+\delta}|s^-|^2A^2-(1+\delta)|ms^+|^2(A+B)^2\\
 &=  \frac{\delta}{1+\delta}|s^-|^2A^2-\delta |ms^+|^2(A+B)^2
 \end{align*}
 for any $\delta>0$. Choosing for instance $\delta:= \frac{2g^3}{\mu_0^+ f^3}$, then
 \begin{align*}
 \frac{1}{12 ms^+}f^3\left(s^- A  + ms^+ (A+B)\right)^2 + \frac{s^+}{3\mu_0^+}g^3(A+B)^2\geq  \frac{s^+}{6\mu_0^+}g^3(A+B)^2+\frac{|s^-|^2f^3g^3}{6 m s^+(\mu_0^+f^3+2g^3)}A^2.
 \end{align*}
 and the assertion follows for
 \[
 \eps= \min\left( \frac{s^+}{6\mu_0^+}g^3,\frac{|s^-|^2f^3g^3}{6 m s^+(\mu_0^+f^3+2g^3)}\right).
 \]
\end{proof}

\section*{Acknowledgement} 
Christina Lienstromberg has been supported by the Deutsche Forschungsgemeinschaft (DFG, German Research Foundation) through the collaborative research centre 'The mathematics of emerging effects' (CRC 1060, Project-ID  211504053) and the Hausdorff Center for Mathematics (GZ 2047/1, Project-ID 390685813). Gabriele Bruell gratefully acknowledges financial support by the Deutsche Forschungsgemeinschaft (DFG) through CRC 1173.
Moreover, the authors are grateful  for fruitful discussions on the topic with Jonas Jansen.
\bigskip

%=============================================================================
%=============================================================================
%=============================================================================
\bigskip

\noindent\textsc{MSC (2010): 76A05, 76A20, 35B40, 35Q35, 35K41, 35K65.}

\medskip

\noindent\textsc{Keywords:} Non-Newtonian fluid, Ellis law, thin-film equation, two-phase flow, degenerate parabolic system, strong solution, local existence and uniqueness, long-time asymptotics
\bigskip


\begin{thebibliography}{10}
%\bibitem{BeFr90}
%F. Bernis and A. Friedman: Higher order nonlinear degenerate parabolic equations, \textit{J. Differential Equations} {\bf{83}} (1990), 179--206.
%
\bibitem{A:1995}
H.~{Amann}.
\newblock {\it Linear and Quasilinear Parabolic Problems, Volume I: Abstract Linear Theory.}
\newblock Birkh\"auser, Basel, 1995.
%
\bibitem{A:1993}
H.~{Amann}.
\newblock {Nonhomogeneous Linear and Quasilinear Elliptic and Parabolic Boundary Value Problems.} In:
\newblock{\it Function Spaces, Differential Operators and Nonlinear Analysis},
\newblock edited by H.-J. Schmeisser, H. Triebel.
\newblock Teubner-Texte zur Math. \textbf{133}, 9--126, Stuttgart, Leipzig 1993.

\bibitem{AG:2002}
L.~{Ansini} and L.~Giacomelli.
\newblock {Shear-thinning liquid films: macroscopic and asymptotic behavior
by quasi-self-similar solutions}.
\newblock {\it Nonlinearity}, \textbf{256}, 2147--2164, 2002.
%
\bibitem{AG:2004}
L.~{Ansini} and L.~Giacomelli.
\newblock {Doubly nonlinear thin-film equations in one space dimension}.
\newblock {\textit Arch. Rational Mech. Anal.}, \textbf{173}, 89--131, 2004.


\bibitem{BBD}
E. Beretta, M. Bertsch and R. Dal Passo. Nonnegative solutions of a fourth-order nonlinear degenerate parabolic equation. \textit{Arch. Rational Mech. Anal. Math. Modeling and Num. Analysis}, \textbf{129}, 175--200, 1995.
%
%

\bibitem{BF}
F. Bernis and A. Friedman. Higher order nonlinear degenerate parabolic equations. \textit{J. Differential Equations}, \textbf{83}(1):179--206, 1990.
%
\bibitem{BP}
A. L. Bertozzi and M. Pugh. The lubrication approximation in thin viscous films: Regularity and long time behavior of weak
solutions. \textit{Comm. Pure Appl. Math}, \textbf{49}, 85--123, 1996.
%
\bibitem{BG}
G. Bruell and R. Granero-Belinch\'{o}n. On the thin film Muskat and the thin film Stokes equations. \textit{J. Math. Fluid Mech.} \textbf{21}, 2019.
%
\bibitem{CC}
E.T. Castellana, P.S. Cremer. Solid supported lipid bilayers: From biophysical studies to sensor design. \textit{Surface Science Reports}  \textbf{61}, 429--444, 2006.
%
\bibitem{DiB} E. DiBenedetto. \textit{Degenerate Parabolic Equations}. Springer-Verlag, New York, 1993.
%
\bibitem{E:1969}
S.~D.{Eidel'man}.
\newblock {\textit{Parabolic Systems.}}
\newblock North-Holland Publishing Company, Amsterdam and Wolters-Noordhoff Publishing, Groningen, 1969.
%
\bibitem{E:2010} L.C.~Evans. 
\newblock{\it Partial Differential Equations.}
\newblock{American Math. Soc., Providence RI, 2010.}
%
\bibitem{EMM}
J. Escher, A.-V. Matioc, B.-V.Matioc. Thin-film approximations of the two-phase Stokes problem. \textit{Nonlinear Anal. Theory Methods Appl.} \textbf{76}, 1--13, (2013)
%
\bibitem{EM}
J. Escher and  B.-V. Matioc. Non-negative global weak solutions for a degenerated parabolic system approximating the two-phase Stokes problem. \textit{J. Differential Equations} \textbf{256}(8), 2659--2676, (2014)
%
\bibitem{Espinosa}
G. Espinosaa, I. L\'{o}pez--Monterob, F. Monroya, and D. Langevina. Shear rheology of lipid monolayers and insights on membrane fluidity. \textit{PNAS} \textbf{108}(15), 2011.
%
\bibitem{GGKO:2014}
L.~Giacomelli, M.~V.~Gnann, H.~Kn\"upfer \& F.~{Otto}.
\newblock{Well-posedness for the Navier-slip thin-film equation in the case of complete wetting.}
\newblock{\it J. Differential Equations}, \textbf{257}, 15--81, 2014.
%
\bibitem{GK:2010}
L.~Giacomelli \& H.~Kn\"upfer.
\newblock{A free boundary problem of fourth order: classical solutions in weighted H\"older spaces.}
\newblock{\it Comm. Partial Differential Equations}, \textbf{35}, 2059--2091, 2010.
%
\bibitem{GKO:2008}
L.~Giacomelli, H.~Kn\"upfer \& F.~{Otto}.
\newblock{Smooth zero-contact-angle solutions to a thin-film equation around the steady state.}
\newblock{\it J. Differential Equations}, \textbf{245}, 1454--1506, 2008.
%
\bibitem{GO}
L. Giacomelli and F. Otto.
\newblock {Rigorous lubrication approximation}.
\newblock {\it Interfaces Free Bound.} \textbf{5}(4):483--529, 2003.
%
\bibitem{GnPet:2018}
M.~V.~Gnann \& M.~Petrache.
\newblock{The Navier-slip thin-film equation for 3D fluid films: existence and uniqueness.}
\newblock{\it J. Differential Equations.}, \textbf{265}, 5832--5958, 2018.
%
\bibitem{GP}
M. G\"unther and G. Prokert.
\newblock {A justification for the thin film approximation of Stokes flow with surface tension}.
\newblock {\it J. Differential Equations.} \textbf{245}(10):2802--2845, 2008.
%
\bibitem{Hermans}
E. Hermans and J. Vermant. Interfacial shear rheology of DPPC under physiologically relevant conditions. \textit{Soft Matter} \textbf{10}, 175--186, 2014.
%
\bibitem{HI}
M. Haragus and G. Iooss.
\newblock {\it Local bifurcation, center manifolds, and normal forms in infinite-dimensional dynamical systems}.
\newblock {Springer-Verlag London, Ltd., London; EDP Sciences, Les Ulis, 2011. } 
%
\bibitem{K1}
J. R. King. The spreading of power-law fluids. \textit{IUTAM Symposium on Free Surface Flows}, 153--160, 2001.
%
\bibitem{K2}
 J. R. King. Two generalisations of the thin film equation. \textit{Math. Comput. Model.}, \textbf{34}:737--756, 2001.
%
\bibitem{Kn:2011}
H.~Kn\"upfer.
\newblock{Well-posedness for the Navier-slip thin-film equation in the case of partial wetting.}
\newblock{\it Arch. Ration. Mech. Anal.}, \textbf{218}, 1083--1130, 2015.
%
\bibitem{Kn:2015}
H.~Kn\"upfer.
\newblock{Well-posedness for a class of thin-film equations with general mobility in the regime of partial wetting.}
\newblock{\it Comm. Pure Appl. Math.}, \textbf{64}, 1263--1296, 2011.
%
\bibitem{LM:2020}
C.~{Lienstromberg} and S.~M\"uller.
\newblock {Local strong solutions to a quasilinear degenerate fourth-order thin-film equation}.
\newblock {\it Nonlinear Differ. Equ. Appl.} \textbf{27}, 16, 2020.
%
\bibitem{L:1995}
A.~Lunardi.
\newblock {\it Ananlytic Semigroups and Optimal Regularity in Parabolic Problems}.
\newblock {Progress in Nonlinear Differential Equations and Their Applications.}, \textbf{16},
\newblock{Birkh\"auser, Basel, 1995.}
%
% \bibitem{MRZ}
% J. M\'{a}lek, K.R. Rajagopal, and J. \v{Z}abensk\'{y}. On power-law fluids with the power-law index proportional to the pressure. \textit{Applied Mathematics Letters}, \textbf{62}, 118--123, 2016.
%
\bibitem{MP}
B.-V. Matioc and G. Prokert.
\newblock{
Hele--Shaw flow in thin threads: a rigorous limit result.}
\newblock{\it Interfaces Free Bound.} \textbf{14}(2), 205–230, 2012. 
%
\bibitem{Milke} A. Mielke. 
\newblock{Reduction of quasilinear alliptic equations in cylindrical domains with applications.}
\newblock{\it Math. Methods Appl. Sci.}, \textbf{10}(1):51--66, 1988.
%
\bibitem{Myres} T.G. Myres. 
\newblock{Application of non-Newtonian models to thin film flow.}
\newblock{\it Physical Review E 72}, 066302, 2005.
%
%
%\bibitem{Tanner}
%R.I. Tanner, K. Walters.
%\newblock {\textit Rheology: a historical perspective}.
%\newblock {Elsevier , Amsterdam, 1998.}
%%
\bibitem{RW}
K.J. Ruschak and S.J. Weinstein. \textit{Coating flows}. Annu. Rev. Fluid Mech. \textbf{36}, 29–53, 2004.
%
\bibitem{Siontorou}
Ch.G. Siontorou, G.-P. Nikoleli, D.P. Nikolelis, and S.K. Karapetis. Artificial Lipid Membranes: Past, Present, and Future. \textit{Membranes} \textbf{7} (38), 2017.	
%
\bibitem{SW}
J.R. SooHoo and G.M. Walker. \textit{Microfluidic aqueous two phase system for leukocyte concentration from whole blood}. Biomed. Microdevices, \textbf{11}:323--329, 2009.

%
\bibitem{WS:1994}
D. E. Weidner and L. W. Schwartz. Contact-line motion of shear-thinning liquids. \textit{Physics of Fluids}, \textbf{6}:3535--3538, 1994.
%

\end{thebibliography}
\end{document}